\documentclass[12 pt]{amsart}
\usepackage{graphicx}

\usepackage{CJKutf8}
\usepackage[dvipsnames]{xcolor}
\usepackage{url}

\usepackage{tikz,tikz-cd,varwidth, amssymb}
\usepackage[hidelinks]{hyperref}
\usetikzlibrary{calc, decorations}
\usepackage[text={6 in, 7.5 in},centering]{geometry}
\usepackage{enumerate, enumitem, stmaryrd}
\usepackage[all,cmtip]{xy}

\usepackage{mathdots}
\usepackage{etoolbox}

\makeatletter
\patchcmd{\@maketitle}{\global\topskip42\p@}{\global\topskip12\p@}{}{}
\patchcmd{\@maketitle}{\dimen@34\p@}{\dimen@22\p@}{}{}
\patchcmd{\@starttoc}{32\p@\@plus14\p@}{16\p@\@plus4\p@}{}{}
\makeatother

\newcommand{\fH}{\mathfrak H}
\newcommand{\fU}{\mathfrak U}

\newcommand{\CC}{\mathbb{C}}

\newcommand{\LL}{\mathbb{L}}

\newcommand{\QQ}{\mathbb{Q}}
\newcommand{\RR}{\mathbb{R}}

\newcommand{\ZZ}{\mathbb{Z}}

\newcommand{\cA}{\mathcal{A}}

\newcommand{\cC}{\mathcal{C}}

\newcommand{\cH}{\mathcal{H}}

\newcommand{\cM}{\mathcal{M}}

\newcommand{\cP}{\mathcal{P}}
\newcommand{\cQ}{\mathcal{Q}}

\newcommand{\cU}{\mathcal{U}}

\newcommand{\BS}{\mathrm{BS}}

\newcommand{\EQ}{{}^{Q}\!E}
\newcommand{\EX}{{}^{X}\!E}

\newcommand{\Set}{\mathsf{Set}}
\newcommand{\op}{\mathrm{op}}

\newcommand{\Gr}{\operatorname{Gr}}
\newcommand{\rank}{\operatorname{rank}}
\newcommand{\Sp}{\mathrm{Sp}}
\newcommand{\Sym}{\operatorname{Sym}}

\newcommand{\St}{\operatorname{St}}
\newcommand{\Cone}{\operatorname{Cone}}
\newcommand{\Tot}{\operatorname{Tot}}

\newcommand{\es}{\operatorname{es}}

\newcommand{\GL}{\mathrm{GL}}
\newcommand{\SL}{\mathrm{SL}}

\newcommand{\rt}{\mathrm{rt}}
\newcommand{\trop}{\mathrm{trop}}

\renewcommand{\>}{\rangle}

\newcommand{\down}[2]{\xymatrix@R=6mm@C=2mm{
#1\ar[d]\\ #2
}}

\newcommand{\downlabel}[3]{\xymatrix@R=6mm@C=2mm{
{#1}\ar[d]^<<<{#3} \\ #2
}}

\newcommand{\squarediagram}[4]{\xymatrix@R=8mm@C=8mm{
#1\ar[d]\ar[r] & #2\ar[d] \\ #3\ar[r] &#4
}}
\newcommand{\squarediagrammapsto}[4]{\xymatrix@R=8mm@C=8mm{
#1\ar@{|->}[d]\ar@{|->}[r] & #2\ar@{|->}[d] \\ #3\ar@{|->}[r] &#4
}}
\newcommand{\squarediagramlabel}[8]{\xymatrix@R=8mm@C=8mm{
#1\ar[d]_{#6}\ar[r]^{#5} & #2\ar[d]^{#7} \\ #3\ar[r]^{#8} &#4
}}

\newcommand{\isocelesdown}[3]{\xymatrix@R=6mm@C=0mm{
& {#1}\ar[dl] \ar[dr] & \\
{#2} \ar[rr] && {#3}
}}

\newcommand{\isocelesdownlabel}[6]{\xymatrix@R=6mm@C=0mm{
& {#1}\ar[dl]_<<<<{#4} \ar[dr]^<<<<{#5} & \\
{#2} \ar[rr]_{#6} && {#3}
}}

\newcommand{\isocelesup}[3]{\xymatrix@R=6mm@C=0mm{
 #1\ar[rr]\ar[dr]  && #2\ar[dl] \\
 & #3 &
}}

\newcommand{\isocelesuplabel}[6]{\xymatrix@R=6mm@C=0mm{
 #1\ar[rr]^{{#4}} \ar[dr]_<<<{#5} && #2\ar[dl]^<<<{#6} \\
 & #3 &
}}

\newtheorem{theorem}{Theorem}[section]
\newtheorem{proposition}[theorem]{Proposition}
\newtheorem{lemma}[theorem]{Lemma}
\newtheorem{conjecture}[theorem]{Conjecture}
\newtheorem{corollary}[theorem]{Corollary}

\newtheorem{definition}[theorem]{Definition}

\theoremstyle{definition}

\newtheorem{remark}[theorem]{Remark}

\newcommand{\id}{\operatorname{id}}

\title[A Hopf algebra action on cohomology]{A Hopf algebra action on the cohomology of moduli of abelian varieties}
\author{Urshita Pal}
\author{Sam Payne}
\date{\today}

\begin{document}

\begin{abstract}
We show that the rational compactly supported cohomology of the moduli space $\cA = \bigsqcup_{g \geq 0} \cA_g$ of principally polarized abelian varieties is free as a bigraded module over the Hopf algebra formed by its weight-zero subspace. The action arises from a symplectic analog of the Waldhausen $S$-construction. This structure yields new lower bounds for Betti numbers and the dimensions of the graded pieces of the weight filtration on the cohomology of $\cA_g$. In particular, $\dim H_c^{2g + k}(\cA_g)$ grows at least exponentially with $g$ for every nonnegative $k \neq 1, 3$. We also give conjectural descriptions of the full mixed Hodge structures in genera four and five, and of mixed Tate subquotients of $H_c^{2g}(\cA_g)$ associated with shuffle algebras and multiple zeta values in every genus.
\end{abstract}

\maketitle

\setcounter{tocdepth}{1}
{\small\tableofcontents}

\section{Introduction}\label{sec:intro}

Let $\cA = \bigsqcup_{g \geq 0} \cA_g$ be the moduli space of principally polarized abelian varieties. It is a smooth complex Deligne--Mumford stack with one component $\cA_g$ for each genus $g \geq 0$. Its compactly supported rational cohomology
\[
H_c^*(\cA) = \bigoplus_{g, k \geq 0} H_c^k(\cA_g)
\]
is bigraded by genus and cohomological degree, with $H_c^k(\cA_g)$ in bidegree $(g, k)$. Each summand carries Deligne's mixed Hodge structure, whose weight filtration is supported in nonnegative weights.

The proper maps $\cA_g \times \cA_h \to \cA_{g + h}$ given by products of principally polarized abelian varieties induce a coproduct on $H_c^*(\cA)$ respecting both gradings. With this coproduct, the weight-zero subspace $W_0 H_c^*(\cA)$ is a bigraded cocommutative Hopf algebra \cite{BCGP}. Our main result extends its left regular action to the full compactly supported cohomology.

\begin{theorem}\label{thm:free-module}
There is a bigraded left action of the Hopf algebra $W_0 H_c^*(\cA)$ on $H_c^*(\cA)$, making $H_c^*(\cA)$ a free bigraded module. The action map is a morphism of coalgebras
\[
\mu \colon W_0 H_c^*(\cA) \otimes H_c^*(\cA) \longrightarrow H_c^*(\cA).
\]
\end{theorem}

\noindent In particular, $H_c^*(\cA)$ is a direct sum of copies of $W_0 H_c^*(\cA)$, shifted in genus and cohomological degree. Theorem~\ref{thm:iota-weight-zero} identifies the action on $1 \in H_c^0(\cA_0)$ with the natural inclusion of the weight-zero subspace.

The Hopf algebra $W_0 H_c^*(\cA)$ contains a tensor sub-Hopf algebra over which it is free, namely the universal enveloping algebra of the motivic Lie algebra of mixed Tate motives over $\ZZ$. Its genus-$g$ part lies in degree $2g$ and has dimension asymptotic to $C \beta^g$, where $C > 0$ and $\beta = 1.3247 \ldots$ is the real root of $t^3 - t - 1$. See Section~\ref{sec:betti-exponential-growth}. Together with Theorem~\ref{thm:free-module}, this gives the following.

\begin{corollary}\label{cor:betti-growth}
For every integer $k$, either $H_c^{2g + k}(\cA_g) = 0$ for all $g$, or
\[
\liminf_{g \to \infty}\left(\dim H_c^{2g + k}(\cA_g)\right)^{1/g} \geq \beta.
\]
The latter is the case for every nonnegative $k \neq 1, 3$.
\end{corollary}

\noindent By Poincar\'e duality, the same statements hold for $H^{g(g - 1) - k}(\cA_g)$.

\begin{remark}
The second part of Corollary~\ref{cor:betti-growth} improves on the growth results obtained from weight-zero cohomology in \cite{BCGP} by using its action on classes in higher degrees, including those detected by pairing with tautological classes. The first part gives a sharp dichotomy between vanishing and exponential growth. Vanishing for $k < 0$ is the symplectic analog of the Church--Farb--Putman conjecture, posed as a question by Br\"uck--Patzt--Sroka \cite[Question~1.1]{BPS}. Our results show that a negative answer would force exponential growth coming from a shifted copy of the universal enveloping algebra of the motivic Lie algebra, whereas the disproof of the Church--Farb--Putman conjecture for mapping class groups was driven by a copy of the motivic Lie algebra itself in $\bigoplus_g W_0 H_c^{2g}(\cM_g)$ \cite{CGP}.
\end{remark}

\bigskip

We construct the action $\mu$ as follows. Recall that the $\GL$-Steinberg homology
\[
\cH := \bigoplus_{n, k} H_k(\GL_n(\ZZ); \St_n)
\]
has a Hopf algebra structure constructed in \cite{AshMillerPatzt24, BCGP}; the two constructions are different, but the Hopf algebras coincide \cite[Theorem~1.1]{PalPayneComparison26}. Here $\St_n$ is the reduced rational homology of the Tits building for $\GL_n$, whose simplices are flags of nonzero proper subspaces of $\QQ^n$. We follow the approach of \cite{BCGP}, in which the Hopf algebra is the $E^1$-page of the homology spectral sequence for the rank filtration on the Waldhausen $S$-construction model of $BK(\ZZ)$. The coproduct is defined via simplicial methods, using a filtered replacement for the diagonal on this model.

Our work here starts similarly, using symplectic variants of both the Waldhausen $S$-construction and Steinberg homology. The symplectic Steinberg homology is
\[
\cH^\omega := \bigoplus_{n, k} H_k(\Sp_{2n}(\ZZ); \St_n^\omega),
\]
where $\St_n^\omega$ is the reduced rational homology of the Tits building of nonzero isotropic subspaces of $\QQ^{2n}$, with its standard symplectic form; we set $\St_0^\omega = \QQ$. We grade $\cH$ by rank and total degree and $\cH^\omega$ by genus and total degree, placing each $H_k$-summand in bidegree $(n, n + k)$; its spectral-sequence bidegree is $(n, k)$.

The orbifold classifying space of $\cA_g$ has homotopy type $B \Sp_{2g}(\ZZ)$. Poincar\'e duality and Borel--Serre duality identify $\cH^\omega$ with the bigraded dual of $H_c^*(\cA)$ through
\[
H_k\bigl(\Sp_{2g}(\ZZ); \St_g^\omega\bigr) \cong H_c^{g + k}(\cA_g)^\vee.
\]

\begin{theorem}\label{thm:main}
The symplectic Steinberg homology $\cH^\omega$ is a bigraded comodule algebra over the $\GL$-Steinberg Hopf algebra $\cH$.
\end{theorem}

\noindent Miller, Pal, and Patzt also construct a comodule-algebra structure on symplectic Steinberg homology for groups over Dedekind domains \cite[Theorem~1.1]{MPP}. Their construction works directly with Steinberg modules and yields lower bounds for the Betti numbers of congruence subgroups of $\Sp_{2n}(\ZZ)$ and of symplectic groups over number rings that are not principal ideal domains.

The product on $\cH^\omega$ is straightforward, induced by block sums of symplectic matrices and joins of Tits buildings. The coaction of $\cH$ is more subtle. We define a symplectic analog of the $S$-construction, denoted $X$, and use simplicial methods to give a filtered map $\rho_X \colon X \to BK(\ZZ) \times X$. The coaction
\[
\rho \colon \cH^\omega \to \cH \otimes \cH^\omega
\]
is induced by $\rho_X$ on the $E^1$-pages of the homology spectral sequences.

Let $e \in H_0(\GL_1(\ZZ); \St_1)$ be the distinguished primitive and set $\overline{\cH} = \cH/(e)$. By \cite{BCGP}, $\overline{\cH}^{\vee} \cong W_0 H_c^*(\cA)$. Composing $\rho$ with the quotient on the first factor gives a coaction $\bar \rho$, whose bigraded dual is the action $\mu$ in Theorem~\ref{thm:free-module}. Dualizing $\rho$ itself makes $H_c^*(\cA)$ a module coalgebra over $\cH^\vee$, but it is easy to see (e.g., for degree reasons) that this larger action is not free.

By the structure theorem for module coalgebras of Milnor and Moore \cite[Theorem~4.4]{MilnorMoore}, freeness reduces to the injectivity of the map
\[
\iota \colon W_0H^*_c(\cA) \to H^*_c(\cA)
\]
given by $a \mapsto \mu(a \otimes 1)$; see Proposition~\ref{prop:MM-freeness}. We prove this by comparing the action with tropical pullback on an $S$-construction analog of a toroidal compactification. A further comparison with the Borel--Serre compactification gives the following result.

\begin{theorem}\label{thm:iota-weight-zero}
The map $\iota$ is the inclusion of the weight-zero subspace in Deligne's mixed Hodge structure on $H_c^*(\cA)$, and the action $\mu$ extends the Hopf algebra product on $W_0 H_c^*(\cA)$.
\end{theorem}

The compatibility between products and coproducts makes Hopf structures powerful tools for proving nonvanishing results and cohomological lower bounds. A recent example is the work of Kupers, Miller, and Patzt \cite[Theorem~A and Corollary~C]{KupersMillerPatzt26}, who use a determinant-twisted variant of $\cH$ to prove that all Morita classes in $H_{4k}(\operatorname{Out}(F_{2k + 2}); \QQ)$, $k \geq 1$, are nonzero. This confirms Morita's conjecture and disproves Kontsevich's finite-dimensionality conjecture for even Lie graph homology.

In a similar spirit, we apply the module coalgebra structure to obtain lower bounds for Betti numbers and the dimensions of the weight-graded pieces of $H_c^*(\cA)$, and to motivate conjectures about its mixed Hodge structure. The coproduct is induced by pullback along product morphisms (Proposition~\ref{prop:block-sums}) and hence is a morphism of mixed Hodge structures. Together with Theorems~\ref{thm:main} and \ref{thm:iota-weight-zero}, this gives the following.

\begin{proposition}\label{prop:weight-convolution}
For all $g, m, w \geq 0$, we have
\[
\dim \Gr_w^W H_c^m(\cA_g) \geq \sum_{h = 0}^g \dim W_0 H_c^{m - w}(\cA_h) \, \dim \Gr_w^W H_c^w(\cA_{g - h}).
\]
\end{proposition}

\noindent These bounds give tantalizing evidence for mixed Tate structures associated with multiple zeta values in $H_c^{2g}(\cA_g)$. We formulate a precise conjecture as follows.

Write $V(r) = V \otimes \QQ(r)$ for the Tate twist of a mixed Hodge structure $V$, where $\QQ(r)$ has Hodge type $(-r, -r)$ and weight $-2r$. Following Brown \cite[\S3.1]{BrownDecomposition}, let $\cU'$ be the graded shuffle Hopf algebra on $f_3, f_5, f_7, \ldots$, with $\deg f_{2i + 1} = 2i + 1$ and deconcatenation coproduct. Regard $\cU'$ as a comodule over itself via its coproduct, and restrict to the elements of degree at most $g$. Write $\cU'_{\leq g}$ for the Hodge realization of the corresponding mixed Tate motive over $\ZZ$. The class $\lambda_1 \cdots \lambda_{g - 1}$ spanning the socle of the tautological ring defines, by Poincar\'e duality, a tautological projection
\[
p_g \colon H_c^{2g}(\cA_g) \longrightarrow \QQ(-g).
\]

\begin{conjecture}\label{conj:multizeta-socle}
For every $g \geq 1$, the mixed Hodge structure $(\cU'_{\leq g})^\vee(-g)$ is a subquotient of $H_c^{2g}(\cA_g)$, with highest-weight quotient $\QQ(-g)$ induced by $p_g$.
\end{conjecture}

\noindent We prove the corresponding lower bounds on the dimensions of the weight graded pieces in every genus. Namely, if $a_n$ is the number of words of total degree $n$ in $f_3, f_5, f_7, \ldots$, then
\[
\dim \Gr_{2(g - n)}^W H_c^{2g}(\cA_g) \geq a_n \qquad (0 \leq n \leq g).
\]
This follows from Proposition~\ref{prop:weight-convolution} and the structure of the tautological projection; see Proposition~\ref{prop:multizeta-weight-bound}. We prove the conjecture for $g \leq 4$.

By Poincar\'e duality, the conjecture is equivalent to $\cU'_{\leq g}(-g(g - 1)/2)$ being a subquotient of $H^{g(g - 1)}(\cA_g)$, with the image of $\lambda_1 \cdots \lambda_{g - 1}$ spanning its lowest-weight subspace. For example, $\zeta(3, 5)$ should occur as a period of a Tate twist of a subquotient of $H^{56}(\cA_8)$. We further expect that, as $g$ varies, the $\QQ[(2 \pi i)^{\pm 1}]$-span of the periods of $H^{g(g - 1)}(\cA_g)$ contains every multiple zeta value and, more generally, every period of a mixed Tate motive over $\ZZ$ \cite[Corollary~1.2]{BrownMixedTate}.

In genera four and five, the bounds lead to predictions for the full mixed Hodge structures. Let $E_3$ denote the mixed Hodge structure in the nonsplit extension
\begin{equation}\label{eq:zeta-three-extension}
0 \longrightarrow \QQ \longrightarrow E_3 \longrightarrow \QQ(-3) \longrightarrow 0
\end{equation}
given by a nonzero rational multiple of $\zeta(3)$.

\begin{conjecture}\label{conj:genus-four-mhs}
The mixed Hodge structures on the compactly supported cohomology groups of $\cA_4$ are
\[
H_c^k(\cA_4) \cong
\begin{cases}
E_3(-1), & k = 8,\\
\QQ(-7)^{ \oplus 2}, & k = 14,\\
\QQ(-k/2), & k \in \{10, 12, 16, 18, 20\},\\
0, & \text{otherwise}.
\end{cases}
\]
\end{conjecture}

\noindent The predicted Betti numbers agree with those suggested by Hulek and Tommasi \cite{HulekTommasiA4}. The degree-$8$ assertion is a consequence of their dimension calculation, the identification $H_c^6(\cA_3) \cong E_3$ in \cite[Theorem~1.19]{BCGP}, and a product-cycle calculation of Canning, Oprea, and Pandharipande \cite{COP}. These show that the coproduct
\[
H_c^8(\cA_4) \longrightarrow H_c^6(\cA_3) \otimes H_c^2(\cA_1)
\]
is an isomorphism; see Proposition~\ref{prop:genus-four-extension}. Together with the known cohomology in genera at most three \cite{HainA3,BCGP}, this proves Conjecture~\ref{conj:multizeta-socle} for $g \leq 4$. The full genus-four conjecture is equivalent to the vanishing of odd cohomology.

For genus five, let $E_{3, 5}$ denote the mixed Hodge structure in an extension
\[
0 \longrightarrow \QQ \longrightarrow E_{3, 5} \longrightarrow \QQ(-3) \oplus \QQ(-5) \longrightarrow 0
\]
whose pullbacks to the two summands are given by nonzero rational multiples of $\zeta(3)$ and $\zeta(5)$, respectively.

\begin{conjecture}\label{conj:genus-five-mhs}
The mixed Hodge structures on the compactly supported cohomology groups of $\cA_5$ are
\[
H_c^k(\cA_5) \cong
\begin{cases}
E_{3, 5}^\vee(-5), & k = 10,\\
E_3(-3), & k = 12,\\
\QQ, & k = 15,\\
\QQ(-k/2), & k \in \{14, 26, 28, 30\},\\
\QQ(-k/2)^{ \oplus 2}, & k \in \{16, 18, 20, 22, 24\},\\
0, & \text{otherwise}.
\end{cases}
\]
\end{conjecture}

Note that Conjectures~\ref{conj:genus-four-mhs} and \ref{conj:genus-five-mhs} are compatible with the $\ell$-adic Euler characteristics computed by Ta\"ibi \cite[Theorem~1]{TaibiEuler}, and the conjectured dimensions of the weight-graded pieces are equal to lower bounds that follow easily from Proposition~\ref{prop:weight-convolution}.

\medskip

\enlargethispage{\baselineskip}
\noindent{\textbf{Acknowledgments.}} SP was supported in part by NSF grant DMS-2542134 and a Simons Fellowship. He carried out part of this work at the Institute for Advanced Study during a visit supported by the Charles Simonyi Endowment.

\section{The \texorpdfstring{$S$}{S}-construction and its symplectic analog}\label{sec:S-constructions}

In this section, we briefly recall the model of $BK(\ZZ)$ produced by the Waldhausen $S$-construction and the Hopf algebra structure on the Quillen spectral sequence of its rank filtration, following \cite[Sections~2--3]{BCGP}. We then present a symplectic analog of this construction with a genus filtration and a filtered coaction of $BK(\ZZ)$.

\subsection{The Waldhausen construction}\label{sec:S}

For each $p \geq 0$, let $S_p$ be the groupoid whose objects are commutative diagrams of finitely generated projective $\ZZ$-modules
\begin{equation} \label{eq:Sdiagram}
\begin{aligned}
\xymatrix {0 \ar[r] & P_{0, 1} \ar[d] \ar[r] & P_{0, 2} \ar[d] \ar[r] & \cdots \ar[d] \ar[r] & P_{0, p} \ar[d]\\
& 0 \ar[r] & P_{1, 2} \ar[d] \ar[r] & \cdots \ar[d] \ar[r]  & P_{1, p} \ar[d] \\
& & 0 \ar[r] & \ddots \ar[d] \ar[r] & \vdots \ar[d] &  \\ & & & 0  \ar[r] & P_{p - 1, p} \ar[d]  \\
& & & & 0,}
\end{aligned}
\end{equation}
in which the top row consists of inclusions with projective cokernels and, for each $0 < i < j \leq p$, the module $P_{i, j}$ is a chosen cokernel of $P_{0, i} \to P_{0, j}$, equipped with its quotient map. The vertical arrows are the resulting quotient maps, and the lower horizontal arrows are induced by the inclusions in the top row. The arrows in $S_p$ are isomorphisms of such diagrams.

Note, in particular, that an object \eqref{eq:Sdiagram} in $S_p$ is determined by its top row, up to choices of cokernels. When no confusion seems possible, we denote such an object more simply by $0 \to P_1 \to \cdots \to P_p$, with $P_i := P_{0, i}$.

We consider functors $d_i \colon S_p \to S_{p - 1}$, for $0 \leq i \leq p$, obtained by deleting the $i$th row and column (i.e., removing the objects labeled $P_{i, j}$ and $P_{k, i}$) and reindexing the remaining terms. Similarly, for $0 \leq i \leq p$, we consider functors $s_i \colon S_p \to S_{p + 1}$ obtained by duplicating the $i$th row and column, inserting identity maps, and reindexing. This makes $S_\bullet$ a simplicial groupoid. Its nerve is a bisimplicial set whose geometric realization is
\[
BK(\ZZ) := |N_\bullet S_\bullet |.
\]
For $n \geq 0$, we filter $BK(\ZZ)$ by rank, defining $F_n S_\bullet$ to be the full simplicial subgroupoid on diagrams whose top-right object $P_p = P_{0, p}$ has rank at most $n$, and setting
\[
F_n BK(\ZZ) := |N_\bullet F_n S_\bullet |.
\]
We set $F_{-1} S_\bullet = \emptyset$ and $F_{-1} BK(\ZZ) = \emptyset$. The \emph{Quillen spectral sequence} $\EQ^*_{*, *}$ is the associated rational homology spectral sequence, converging to $H_*(BK(\ZZ); \QQ)$ and indexed by
\[
\EQ^1_{s, t} := H_{s + t}\bigl(F_s BK(\ZZ), F_{s - 1} BK(\ZZ); \QQ\bigr) \cong H_t(\GL_s(\ZZ); \St_s).
\]
For $s \geq 1$, the Tits building $T_s(\QQ)$ is the nerve of the poset of nonzero proper subspaces of $\QQ^s$, and
\[
\St_s := \widetilde H_{s - 2}\bigl(T_s(\QQ); \QQ\bigr)
\]
is the Steinberg module, with its natural $\GL_s(\ZZ)$-action. We follow the convention that $\widetilde H_{-1}(\emptyset; \QQ) = \QQ$, so $\St_1 = \QQ$.  We also set $\St_0 = \QQ$.

\medskip

The first page $\EQ^1_{*, *}$ is a connected, bigraded Hopf algebra whose product is graded-commutative with respect to total degree \cite[Theorem~3.18]{BCGP}: for $x \in \EQ^1_{s, t}$ and $y \in \EQ^1_{s', t'}$,
\[
xy = (-1)^{(s + t)(s' + t')} yx.
\]
The multiplication and comultiplication are induced by a filtered product and coproduct on $BK(\ZZ)$.

The product $m \colon BK(\ZZ) \times BK(\ZZ) \to BK(\ZZ)$ is given on objects of $S_p$ by term-by-term direct sums of diagrams, i.e.,
\[
\thickmuskip=3mu
\medmuskip=3mu
\bigl((0 \to P_1 \to \cdots \to P_p), (0 \to P'_1 \to \cdots \to P'_p)\bigr) \longmapsto (0 \to P_1 \oplus P'_1 \to \cdots \to P_p \oplus P'_p).
\]
It is filtered because $m(F_a BK(\ZZ) \times F_b BK(\ZZ)) \subseteq F_{a + b} BK(\ZZ)$, and it is commutative up to a filtration-preserving homotopy; see \cite[Propositions~3.13 and~3.16]{BCGP}.

We use the following conventions for filtered spaces. A map $f \colon Z \to Z'$ is \emph{filtered} if $f(F_n Z) \subseteq F_n Z'$ for every $n$. A homotopy $H \colon [0, 1] \times Z \to Z'$ is \emph{filtered} if
\[
H([0, 1] \times F_n Z) \subseteq F_n Z' \qquad \text{for every } n.
\]
Maps joined by such a homotopy induce the same map on $E^1$ of the associated spectral sequences. Products of filtered spaces carry the total filtration
\[
F_n(Z \times Z') := \bigcup_{a + b \leq n} F_a Z \times F_b Z'.
\]

\medskip

The diagonal $\Delta \colon BK(\ZZ) \to BK(\ZZ) \times BK(\ZZ)$ is not filtered. The coproduct is instead given by a filtered map that is homotopic to $\Delta$. The resulting coproducts on the pages of the Quillen spectral sequence are not all cocommutative. On the abutment, however, this map induces the usual cocommutative homology coproduct. The filtered coproduct is defined via edgewise subdivision, as follows.

Consider the bisimplicial set $Y_{p, q} := N_q S_p$. We take edgewise subdivision in the $p$-direction, setting
\[
\es(Y)_{p, q} := Y_{2p + 1, q} = Y([p] \sqcup [p], [q]),
\]
where $[p] \sqcup [p]$ is ordered by concatenation, i.e., lexicographically on $[1] \times [p]$. There is a natural homeomorphism of geometric realizations
\[
| \es(Y)_{\bullet, \bullet} | \cong |Y_{\bullet, \bullet} | = BK(\ZZ),
\]
and a natural transformation
\[
\es(Y)_{p, q} \to Y_{p, q} \times Y_{p, q}
\]
induced by the two order-preserving inclusions $[p] \to [2p + 1]$ onto the first and last $p + 1$ elements, respectively. The resulting map of spaces
\[
BK(\ZZ) \cong | \es(Y)_{\bullet, \bullet} | \to BK(\ZZ) \times BK(\ZZ)
\]
is naturally homotopic to the diagonal \cite[\S3]{BCGP}. Unlike the diagonal, it is filtered: for every $n$, its restriction to $F_n BK(\ZZ)$ has image in
\[
F_n\bigl(BK(\ZZ) \times BK(\ZZ)\bigr) = \bigcup_{\substack {s, t \geq 0\\s + t \leq n}} F_s BK(\ZZ) \times F_t BK(\ZZ).
\]
Indeed, at the level of groupoids, the $p$th component $\es(S_\bullet)_p = S_{2p + 1} \to S_p \times S_p$ is given on objects by
\[
\begin{aligned}
& (0 \to P_1 \to \cdots \to P_{2p + 1})\\
& \qquad \longmapsto \bigl((0 \to P_1 \to \cdots \to P_p), (0 \to P_{p + 2}/P_{p + 1} \to \cdots \to P_{2p + 1}/P_{p + 1})\bigr).
\end{aligned}
\]
The ranks of these modules satisfy
\[
\rank(P_p) + \rank(P_{2p + 1}/P_{p + 1}) \leq \rank(P_{p + 1}) + \rank(P_{2p + 1}/P_{p + 1}) = \rank(P_{2p + 1}),
\]
because $P_p \subseteq P_{p + 1}$. The natural homeomorphism between the realizations before and after edgewise subdivision preserves the induced filtrations \cite[Lemma~3.3]{BCGP}, so this calculation gives the asserted filtered map on realizations.

\subsection{The symplectic construction}\label{sec:symplectic-S}

We now present a symplectic analog of the Waldhausen $S$-construction model for $BK(\ZZ)$, with projective $\ZZ$-modules replaced by projective $\ZZ$-modules with alternating bilinear pairings. We define a simplicial groupoid of isotropic flags in symplectic lattices, filter its realization by genus, and identify the $E^1$-page of the resulting spectral sequence with symplectic Steinberg homology. The resulting filtered space, with its structure as the geometric realization of a filtered bisimplicial set, is essential for the proofs of our main results.

Consider the category $\cC$ of pairs $(Q, \omega)$ where $Q$ is a finitely generated projective $\ZZ$-module and $\omega \colon Q \times Q \to \ZZ$ is an alternating bilinear pairing. A morphism in this category from $(Q, \omega)$ to $(Q', \omega')$ is a map of $\ZZ$-modules $g \colon Q \to Q'$ such that $\omega = g^* \omega'$. When no confusion seems possible, we write $Q$ for the object $(Q, \omega)$. An injection of finitely generated projective $\ZZ$-modules is \emph{saturated} if its cokernel is torsion-free, or equivalently projective. For each $p \geq 0$, let $\cC_p$ be the groupoid whose objects are commutative diagrams in $\cC$
\begin{equation} \label{eq:Qdiagram}
\begin{aligned}
\xymatrix {0 \ar[r] & Q_{0, 1} \ar[d] \ar[r] & Q_{0, 2} \ar[d] \ar[r] & \cdots \ar[d] \ar[r] & Q_{0, p + 1} \ar[d] \ar[r] &  \cdots \ar[d] \ar[r] & Q_{0, 2p} \ar[d] \ar[r] & Q_{0, 2p + 1} \\  & 0 \ar[r] & Q_{1, 2} \ar[d] \ar[r] & \cdots \ar[d] \ar[r]  & Q_{1, p + 1} \ar[d] \ar[r] & \cdots \ar[d] \ar[r] & Q_{1, 2p} & \\
& & 0 \ar[r] & \ddots \ar[d] \ar[r] & \vdots \ar[d] \ar[r] &  \iddots   & &\\ & & & 0  \ar[r] & Q_{p, p + 1}  & & &}
\end{aligned}
\end{equation}
where $Q_{0, 0} = 0$, the pairing on $Q_{0, 2p + 1}$ is perfect, each arrow $Q_{0, i} \to Q_{0, i + 1}$ is a saturated injection, and the images of $Q_{0, i}$ and $Q_{0, 2p + 1 - i}$ are orthogonal complements in $Q_{0, 2p + 1}$ for $0 \leq i \leq p$. For $0 < i < j \leq 2p + 1 - i$, the module $Q_{i, j}$ is a chosen cokernel, in the underlying category of $\ZZ$-modules, of the map $Q_{0, i} \to Q_{0, j}$, with the induced pairing. This pairing is well-defined because $Q_{0, j} \subseteq Q_{0, 2p + 1 - i} = Q_{0, i}^{\perp}$. The vertical arrows are the resulting quotient maps, and the lower horizontal arrows are induced by the inclusions in the top row. It follows that the pairing on $Q_{i, 2p + 1 - i}$ is perfect and alternating, all horizontal arrows are saturated inclusions, and the images of $Q_{i, j}$ and $Q_{i, 2p + 1 - j}$ are orthogonal complements in $Q_{i, 2p + 1 - i}$ for $i < j < 2p + 1 - i$. The arrows in $\cC_p$ are isomorphisms of such diagrams.

An object of $\cC_p$ is therefore determined, up to choices of cokernels, by its ambient symplectic module $Q_{0, 2p + 1}$ and a flag of saturated isotropic submodules
\[
0 = Q_{0, 0} \subseteq \cdots \subseteq Q_{0, p} \subseteq Q_{0, 2p + 1}.
\]
The remaining terms in the top row are the corresponding orthogonal complements.

We define functors $d_i \colon \cC_p \to \cC_{p - 1}$ for $0 \leq i \leq p$ by deleting the $i$th row, the $i$th column, and the $(2p + 1 - i)$th column (i.e., removing all terms of the form $Q_{i, j}$, $Q_{k, i}$, or $Q_{\ell, 2p + 1 - i}$), and then reindexing the remaining terms. Similarly, for $0 \leq i \leq p$, we define functors $s_i \colon \cC_p \to \cC_{p + 1}$ by duplicating the $i$th row, the $i$th column, and the $(2p + 1 - i)$th column, inserting identity morphisms, and then reindexing. These functors make $\cC_\bullet$ a simplicial groupoid. Taking the nerve yields a bisimplicial set
\[
N_\bullet \cC_\bullet \colon \Delta^\op \times \Delta^\op \to \Set
\]
given by $([p], [q]) \mapsto N_q \cC_p$.
\begin{remark}\label{rem:symplectic-small-model}
The category $\cC$ is essentially small but not small. In order for $N_\bullet \cC_\bullet$ to be a bisimplicial \emph{set}, we tacitly replace $\cC$ with a suitable equivalent small symmetric monoidal category, as in \cite[Remark~2.13]{BCGP}.
\end{remark}

We filter $\cC_p$ by the genus of the ambient symplectic module $Q_{0, 2p + 1}$, where the genus of a module with a perfect alternating pairing is half its rank. More precisely, $F_n \cC_\bullet$ is the full simplicial subgroupoid on diagrams for which $\rank(Q_{0, 2p + 1}) \leq 2n$. We set $F_{-1} \cC_\bullet = \emptyset$. Genus is non-increasing under the simplicial operators and is preserved by arrows in $\cC_\bullet$, so this gives an exhaustive filtration of $N_\bullet \cC_\bullet$. Set
\[
X := |N_\bullet \cC_\bullet |, \qquad F_n X := |N_\bullet F_n \cC_\bullet |.
\]

Termwise orthogonal direct sum defines a map of bisimplicial sets and hence a product
\[
m_X \colon X \times X \longrightarrow X.
\]
It is filtered because $m_X(F_n X \times F_{n'} X) \subseteq F_{n + n'} X$, and its unit is the zero diagram in $F_0 X$. The canonical associativity, unit, and symmetry isomorphisms for orthogonal direct sum preserve genus and induce filtered homotopies on realizations. Hence $m_X$ induces a product on the associated spectral sequence.

\subsection{Symplectic Steinberg homology}\label{sec:symplectic-S-homology}

Let $\EX^*_{*, *}$ be the rational homology spectral sequence of this filtration, indexed by
\[
\EX^1_{n, k} := H_{n + k}(F_n X, F_{n - 1} X; \QQ) \Longrightarrow H_{n + k}(X; \QQ).
\]
For $n \geq 1$, let $T^\omega_n(\QQ)$ be the nerve of the poset of nonzero isotropic subspaces of $\QQ^{2n}$ with the standard symplectic form $\omega$, characterized on the standard basis by
\[
\omega(e_i, f_j) = \delta_{ij}, \qquad \omega(e_i, e_j) = \omega(f_i, f_j) = 0.
\]
The reduced rational homology of $T_n^\omega(\QQ)$ is concentrated in degree $n - 1$, and the symplectic Steinberg module is
\[
\St_n^\omega := \widetilde H_{n - 1}\bigl(T_n^\omega(\QQ); \QQ\bigr).
\]
It carries the natural $\Sp_{2n}(\ZZ)$-action. We set $T_0^\omega(\QQ) = \emptyset$ and $\St_0^\omega = \QQ$.

For each $n \geq 0$, write $\widetilde C_*(T_n^\omega(\QQ); \QQ)$ for the augmented simplicial chain complex computing reduced rational homology. 

We follow the sign conventions of \cite[\S\S2.7 and~8.5]{Weibel}. In particular, the differential on a tensor product of chain complexes is given by $$d(x \otimes y) = dx \otimes y + (-1)^{|x|} x \otimes dy.$$ For bisimplicial chains $A_{p, q}$, with $d_h$ and $d_v$ the alternating sums of horizontal and vertical faces, respectively, the total differential is
\[
d_{\Tot} = d_h + (-1)^p d_v \quad \mbox{on } A_{p, q}.
\]
We give $H_k(\GL_n(\ZZ); \St_n)$ and $H_k(\Sp_{2n}(\ZZ); \St_n^\omega)$ total degree $n + k$ by placing the Steinberg coefficient modules in chain degree $n$. In computing derived coinvariants, we place the group resolution first in the tensor product. All cross products and K\"unneth identifications involving these groups use the resulting Koszul signs. In particular, on spectral-sequence bidegrees $(n, k)$ and $(n', k')$, the cross product is $(-1)^{nk'}$ times the ordinary group-homology cross product with coefficients placed in degree zero.

Let $\Sigma C$ denote the suspension of a chain complex $C$, so $(\Sigma C)_j = C_{j - 1}$ with differential $d_{\Sigma C} = -d_C$. The product on $\EX^1$ is described by the shuffle map
\[
\operatorname{EZ}_{n, n'} \colon \Sigma \widetilde C_*(T^\omega_n(\QQ); \QQ) \otimes \Sigma \widetilde C_*(T^\omega_{n'}(\QQ); \QQ) \longrightarrow \Sigma \widetilde C_*(T^\omega_{n + n'}(\QQ); \QQ),
\]
as follows. Let $0 \subsetneq U_1 \subsetneq \cdots \subsetneq U_r$ and $0 \subsetneq W_1 \subsetneq \cdots \subsetneq W_s$ be strict flags of isotropic subspaces of $\QQ^{2n}$ and $\QQ^{2n'}$, respectively, where $r, s \geq 0$, and set $U_0 = W_0 = 0$. Every monotone lattice path from $(0, 0)$ to $(r, s)$, increasing one coordinate by one at each step, determines the strict flag formed by the nonzero subspaces $U_i \oplus W_j$ encountered along the path. We regard these subspaces as lying in $\QQ^{2(n + n')}$ through the standard symplectic identification $\QQ^{2n} \perp \QQ^{2n'} \cong \QQ^{2(n + n')}$ obtained by concatenating the $e$-vectors and then the $f$-vectors. The term corresponding to a path $\gamma$ occurs with sign $(-1)^{\epsilon(\gamma)}$, where $\epsilon(\gamma)$ is the number of pairs in which a vertical step of $\gamma$ precedes a horizontal step, i.e., the number of inversions in the corresponding $(r, s)$-shuffle.

Note that $\Sigma \widetilde C_*(T_n^\omega(\QQ); \QQ)$ is the normalized relative chain complex of the cone on $T_n^\omega(\QQ)$ relative to the building, represented by flags beginning with $0$. Under this identification, $\operatorname{EZ}_{n, n'}$ is the ordinary Eilenberg--Zilber shuffle followed by orthogonal direct sum, and hence is a chain map. See \cite[\S8.5]{Weibel}.

For $n \geq 1$, let $T_n = T_n^\omega(\QQ)$. A rational isotropic flag $U_0 \subseteq \cdots \subseteq U_p$ determines saturated lattices $I_j = U_j \cap \ZZ^{2n}$. These give an object in $\cC_p$ with ambient lattice $I_0^\perp/I_0$ and isotropic flag $0 \subseteq I_1/I_0 \subseteq \cdots \subseteq I_p/I_0$. Let $\Gamma = \Sp_{2n}(\ZZ)$. Note that $\Cone(T_n)$ is the realization of the poset of all rational isotropic subspaces, including $0$. Sending each flag to the diagram above and each $\gamma \in \Gamma$ to the induced isomorphism of diagrams is compatible with faces and degeneracies. In particular, the face $d_0$ replaces the ambient lattice $I_0^\perp/I_0$ by $I_1^\perp/I_1$. Taking nerves and realizations gives a map $\kappa_n \colon E \Gamma \times_\Gamma \Cone(T_n) \to F_n X$. Since $I_0^\perp/I_0$ has genus $n - \dim U_0$, this map carries $E \Gamma \times_\Gamma T_n$ into $F_{n - 1} X$.

\begin{lemma}\label{lem:building-filtration}
For every $n \geq 1$, pullback by $\kappa_n$ gives an isomorphism in rational cohomology
\[
\kappa_n^* \colon H^*(F_n X, F_{n - 1} X) \xrightarrow{\sim} H_{\Gamma}^*(\Cone(T_n), T_n).
\]
\end{lemma}

\begin{proof}
Let $\cC_p^n$ be the full subgroupoid of $F_n \cC_p$ on diagrams whose ambient symplectic module has genus $n$. In each bidegree, the pointed filtration quotient is
\[
F_n N_q \cC_p/F_{n - 1} N_q \cC_p \cong N_q \cC_p^n \cup \{\infty\}.
\]
The face maps $d_i$ for $i > 0$ preserve genus, while $d_0$ does so only if $Q_{0, 1} = 0$.  If $Q_{0,1} \neq 0$ then $d_0$ decreases genus and sends the simplex to the basepoint $\infty$. Write $Z_{p, q} = \QQ\<N_q \cC_p^n\>$ for the reduced bisimplicial chains. The total complex $\Tot Z_{*, *}$ computes the relative homology of $(F_n X, F_{n - 1} X)$. 

Let $\widetilde \cC_p^n$ be the groupoid whose objects are diagrams in $\cC_p^n$ equipped with a framing $\varphi \colon Q_{0, 2p + 1} \xrightarrow{\sim} (\ZZ^{2n}, \omega)$. Morphisms preserve the framings. The group $\Gamma$ acts freely by postcomposition on the framings, and forgetting the framing identifies $N_q \cC_p^n$ with the orbit set of $N_q \widetilde \cC_p^n$. There is a map
\[
N_q \widetilde \cC_p^n \longrightarrow \pi_0(N_\bullet \widetilde \cC_p^n)
\]
taking each simplex to its connected component. For fixed $p$, these maps form a map of simplicial sets as $q$ varies. An automorphism of a framed diagram is the identity on the ambient module, and hence on every term of the diagram. Thus all automorphism groups are trivial, so the assembled map of simplicial sets is a weak homotopy equivalence.

The connected components of $N_\bullet \widetilde \cC_p^n$ correspond, via the framing, to flags of saturated isotropic sublattices $0 = I_0 \subseteq \cdots \subseteq I_p$ in $\mathbb Z^{2n}$. After adjoining the basepoint, these components form a simplicial set whose degeneracies repeat terms of the flag. Its reduced normalized chains are therefore generated by strict flags $0 \subsetneq I_1 \subsetneq \cdots \subsetneq I_p$. On these generators, $d_0$ is zero and $d_i$ deletes $I_i$ for $i > 0$. Taking rational spans identifies these flags with the simplices of the cone on $T_n$ that contain the vertex $0$. Thus the reduced normalized chain complex is the relative chain complex of $(\operatorname{Cone}(T_n), T_n)$, identified with $\Sigma\widetilde C_*(T_n; \mathbb Q)$ with differential $-d$.

The preceding comparisons give a $\Gamma$-equivariant quasi-isomorphism from the total complex of reduced chains of the framed diagrams to $\Sigma \widetilde C_*(T_n; \QQ)$. Since $\Gamma$ acts freely on framed simplices, this total complex consists of free $\mathbb Q[\Gamma]$-modules. Its coinvariants are $\Tot Z_{*, *}$, which therefore computes the equivariant homology of $(\operatorname{Cone}(T_n), T_n)$ as well as the relative homology of $(F_n X, F_{n-1} X)$. Under these comparisons, $\kappa_n$ induces the identity on rational isotropic flags and hence an isomorphism on homology. Dualizing over $\QQ$ proves the lemma.
\end{proof}

\begin{proposition}\label{prop:E1X}
The spectral sequence $\EX^*_{*, *}$ converges to $H_*(X; \QQ)$, and there is a canonical isomorphism
\[
\EX^1_{n, k} \cong H_k(\Sp_{2n}(\ZZ); \St^\omega_n).
\]
Under this isomorphism, the product on $\EX^1$ is the group-homology cross product followed by the map induced by the block-sum homomorphism $\Sp_{2n}(\ZZ) \times \Sp_{2n'}(\ZZ) \to \Sp_{2n + 2n'}(\ZZ)$ and the equivariant map $\St_n^\omega \otimes \St_{n'}^\omega \longrightarrow \St_{n + n'}^\omega$ induced by $\operatorname{EZ}_{n, n'}$.
\end{proposition}

\begin{proof}
For $n \geq 1$, Lemma~\ref{lem:building-filtration} identifies the relative homology of $(F_n X, F_{n - 1} X)$ with the equivariant homology of $(\Cone(T_n), T_n)$. The latter pair has homology concentrated in degree $n$, with coefficient module $\St_n^\omega$, so the equivariant homology spectral sequence gives the asserted isomorphism. The case $n = 0$ follows from the contractibility of $F_0 X$. Exhaustiveness and the first-quadrant condition give convergence.

As in the proof of \cite[Proposition~2.15]{BCGP}, the framed comparison is compatible with direct sums. Orthogonal direct sum induces block sum on the acting symplectic groups and, after applying the Eilenberg--Zilber shuffle, the map $\operatorname{EZ}_{n, n'}$ on suspended building chains. Passing to derived coinvariants gives the asserted product, with the sign $(-1)^{nk'}$ arising from moving the first coefficient factor past the second group-resolution factor.
\end{proof}

\section{The comodule algebra structure}\label{sec:comod}

We retain the notation $X = |N_\bullet \cC_\bullet |$, the genus filtration $F_n X = |N_\bullet F_n \cC_\bullet |$, and the product $m_X \colon X \times X \to X$ from the previous section, and write $\cH^{\omega} \cong \EX^1$ for the symplectic Steinberg homology. In this section we construct the coaction using filtered maps and homotopies, then verify the comodule algebra axioms on $E^1$.

Recall from \cite[\S3]{BCGP} that edgewise subdivision gives, for every bisimplicial set $Z$, a map
\[ \Phi_Z \colon |Z| \xrightarrow{\ \cong \ } | \es(Z) | \longrightarrow |Z| \times |Z|, \]
homotopic to the diagonal and natural in the sense that $\Phi_W \circ |g| = (|g| \times |g|) \circ \Phi_Z$ for every $g \colon Z \to W$. Applying this construction to $N_\bullet \cC_\bullet$ and $N_\bullet S_\bullet$, we write
\[
\Phi_X \colon X \longrightarrow X \times X, \qquad \Phi_{BK} \colon BK(\ZZ) \longrightarrow BK(\ZZ) \times BK(\ZZ)
\]
for the resulting maps of realizations. Let
\[
f \colon N_\bullet \cC_\bullet \longrightarrow N_\bullet S_\bullet
\]
be the forgetful map of bisimplicial sets that records the underlying isotropic flag and its chosen subquotients, discarding the ambient symplectic form and the orthogonal complements, and use the same letter for its realization $f \colon X \to BK(\ZZ)$. The coaction is
\[ \rho_X := (f \times \id_X) \circ \Phi_X \; \colon \; X \longrightarrow BK(\ZZ) \times X. \]
Write $Q_i = Q_{0, i}$ for the top-row terms. The bisimplicial map defining $\rho_X$ sends $0 \to Q_1 \to \cdots \to Q_{4p + 3}$ in $\es(N_\bullet \cC_\bullet)_{p, 0} = N_0 \cC_{2p + 1}$ to
\[ \bigl(0 \to Q_1 \to \cdots \to Q_p\bigr) \ \times \ \bigl(0 \to Q_{p + 2}/Q_{p + 1} \to \cdots \to Q_{3p + 2}/Q_{p + 1}\bigr). \]
Thus $\rho_X$ retains the initial isotropic flag as a flag of projective modules and passes in the second factor to the symplectic reduction $Q_{p + 1}^{\perp}/Q_{p + 1}$.

\begin{lemma}\label{lem:f-strict}
The filtered map of bisimplicial sets $f$ is strictly monoidal:
\[
f(D \perp D') = f(D) \oplus f(D').
\]
Its realization satisfies $\Phi_{BK} \circ f = (f \times f) \circ \Phi_X$.
\end{lemma}

\begin{proof}
On objects, $f$ sends a $\cC_p$-diagram to the $S_p$-diagram consisting of the isotropic flag $0 \subseteq Q_1 \subseteq \cdots \subseteq Q_p$ and the chosen cokernels $Q_{i, j}$ for $0 < i < j \leq p$. For $1 \leq i \leq p$, the face $d_i$ of $\cC_\bullet$ deletes columns $i$ and $2p + 1 - i$. Only column $i$ belongs to the isotropic flag, so $d_i$ acts on this flag as the corresponding face of $S_\bullet$. The face $d_0$ replaces it by the quotient flag with terms $Q_{1, j} = Q_j/Q_1$, again matching $d_0$ of $S_\bullet$. Degeneracies repeat the corresponding terms in both flags. Thus $f$ commutes with the simplicial operators.

The rank of the isotropic flag is at most the genus of its ambient symplectic module, so $f$ preserves the filtrations. It commutes strictly with direct sums, and naturality of $\Phi$ gives the final identity.
\end{proof}

\begin{lemma}\label{lem:rhofilt}
The coaction $\rho_X$ is filtered.
\end{lemma}

\begin{proof}
Let $0 \to Q_1 \to \cdots \to Q_{4p + 3}$ be a diagram as above with $\rank(Q_{4p + 3}) = 2n$. The initial flag $0 \to Q_1 \to \cdots \to Q_p$ has rank $\rank(Q_p)$. Since $Q_{3p + 2} = Q_{p + 1}^{\perp}$, the symplectic quotient $Q_{3p + 2}/Q_{p + 1}$ has genus $n - \rank(Q_{p + 1})$. Thus
\[
\rank(Q_p) + n - \rank(Q_{p + 1}) = n - \rank(Q_{p + 1}/Q_p) \leq n.
\]
Hence $\rho_X$ preserves the filtration.
\end{proof}

\begin{proposition}\label{prop:coassoc}
There is a filtered homotopy joining the two maps
\[ (\Phi_{BK} \times \id_X) \circ \rho_X \qquad \text{and} \qquad (\id_{BK} \times \rho_X) \circ \rho_X \; \colon \; X \longrightarrow BK(\ZZ) \times BK(\ZZ) \times X. \]

\end{proposition}

\begin{proof}
Let $\psi \colon [0, 1] \to [0, 1]$ be the increasing piecewise-linear map that carries the intervals $[0, \tfrac14]$, $[\tfrac14, \tfrac12]$, and $[\tfrac12, 1]$ linearly onto $[0, \tfrac12]$, $[\tfrac12, \tfrac34]$, and $[\tfrac34, 1]$, respectively. Represent the simplex in the $p$-direction as
\[
\Delta^p = \{(t_1, \ldots, t_p) \in [0, 1]^p : t_1 \leq \cdots \leq t_p\}.
\]
On each $\Delta^p \times \Delta^q$, apply $\psi$ coordinatewise to the first factor and the identity to the second. These maps define a natural self-map $\psi_X$ of $X$. As in the proof of \cite[Proposition~3.11]{BCGP}, let
\[
G_X \colon [0, 1] \times X \longrightarrow X
\]
be the simplexwise straight-line homotopy from $\id_X$ to $\psi_X$. This homotopy is filtered. Moreover, $(\id_X \times \Phi_X) \circ \Phi_X \circ G_X$ is a homotopy from $(\id_X \times \Phi_X) \circ \Phi_X$ to $(\Phi_X \times \id_X) \circ \Phi_X$. By the definition of $\rho_X$ and Lemma~\ref{lem:f-strict}, we have
\[ \begin{aligned}
(\id_{BK} \times \rho_X) \circ \rho_X & = (f \times f \times \id_X) \circ (\id_X \times \Phi_X) \circ \Phi_X,\\
(\Phi_{BK} \times \id_X) \circ \rho_X & = (f \times f \times \id_X) \circ (\Phi_X \times \id_X) \circ \Phi_X.
\end{aligned} \]
Then, since $\rho_X$ and $G_X$ are filtered, the composition 
\[
(\id_{BK} \times \rho_X) \circ \rho_X \circ G_X
\]
is a filtered homotopy, and the proposition follows.
\end{proof}

We next verify the counit identity. Let ${\varepsilon \colon BK(\ZZ) \to \ast}$ be the map to a point. Then $(\varepsilon \times \id_X) \circ \rho_X = \pi_2 \circ \Phi_X$. The second component of $\Phi_X$ is naturally homotopic to $\id_X$ through the simplexwise straight-line homotopy of \cite[\S3]{BCGP}; the homotopy is filtered because it stays within each simplex.

For compatibility with multiplication, note that edgewise subdivision commutes with products, and the natural homeomorphism on realizations is compatible with products. Hence, under the canonical identification of the realization of $(N_\bullet \cC_\bullet) \times (N_\bullet \cC_\bullet)$ with $X \times X$, the map $\Phi_{X \times X}$ is $\tau_{23} \circ (\Phi_X \times \Phi_X)$, where
\[
\tau_{23}\bigl((x_1, x_2), (x_3, x_4)\bigr) = \bigl((x_1, x_3), (x_2, x_4)\bigr)
\]
interchanges the second and third factors. Naturality of $\Phi$ applied to the bisimplicial product defining $m_X$, together with the strict monoidality of $f$, therefore gives
\begin{equation} \label{eq:rhom}
\rho_X \circ m_X = (m \times m_X) \circ \tau_{23} \circ (\rho_X \times \rho_X),
\end{equation}
where $m$ is the product on $BK(\ZZ)$.

\begin{proof}[Proof of Theorem~\ref{thm:main}]
The exterior cross product and the rational K\"unneth theorem give a natural isomorphism 
\[
E^1\bigl(BK(\ZZ) \times X\bigr) \cong \EQ^1 \otimes \EX^1,
\]
as in \cite[\S3]{BCGP}. Since $\rho_X$ is filtered, it induces, through this isomorphism and the identifications $\cH \cong \EQ^1$ and $\cH^\omega \cong \EX^1$, a bigraded map
\[
\rho \colon \cH^\omega \longrightarrow \cH \otimes \cH^\omega.
\]
The preceding filtered homotopies give coassociativity and counitality on $E^1$, and equation~\eqref{eq:rhom} gives compatibility with multiplication. Since $\rho_X$ sends the zero diagram to the pair of zero diagrams, $\rho$ also preserves the unit. Thus $\rho$ makes $\cH^\omega$ a bigraded comodule algebra over $\cH$.
\end{proof}

For the applications to cohomology of moduli of abelian varieties, we pass to a Hopf quotient of $\cH$ and dualize. Recall from \cite[\S1]{BCGP} the distinguished primitive $e \in \EQ^1_{1, 0}$, which is in bidegree $(1,1)$ for the bigrading on $\cH$ and generates a Hopf ideal $(e)$. Setting $\overline{\cH} := \cH/(e)$, we have a surjection of Hopf algebras
\[
q \colon \cH \twoheadrightarrow \overline{\cH}.
\]
Let $\bar \rho := (q \otimes \id_{\cH^\omega}) \circ \rho$. Then
\[
\bar \rho \colon \cH^{\omega} \longrightarrow \overline{\cH} \otimes \cH^{\omega}
\]
is a coaction making $\cH^\omega$ a comodule algebra over the Hopf algebra $\overline \cH$.

We have bigraded isomorphisms
\[
\overline{\cH}^{\vee} \cong W_0 H^*_c(\cA) \quad \mbox{ and } \quad (\cH^\omega)^\vee \cong H^*_c(\cA).
\]
The first isomorphism follows from \cite{BCGP}, and the second from Borel--Serre duality \cite[Theorems~11.4.2 and~11.4.4]{BorelSerre} and Poincar\'e duality.

For the isomorphism $\overline{\cH}^{\vee} \cong W_0 H_c^*(\cA)$, we use the comparison of \cite[\S2.4]{BCGP} with the convention specified in Section~\ref{sec:tropical-comparison}; see \cite[Remark~2.21]{BCGP}.

The product on $\cH^\omega$ dualizes to a coproduct on $H_c^*(\cA)$. Proposition~\ref{prop:block-sums} identifies this coproduct with pullback along the proper product morphism $\cA \times \cA \to \cA$, which shows that it is a morphism of mixed Hodge structures.

\medskip

A graded left \emph{module coalgebra} over a Hopf algebra $A$ is a graded left $A$-module $B$ with $A$-linear coproduct and counit, where $A$ acts on $B \otimes B$ through its coproduct. Here connected means nonnegatively graded with degree-zero part equal to the ground field.

\begin{corollary}\label{cor:module-coalgebra}
The bigraded dual
\[
\mu := \bar \rho^\vee \colon W_0 H^*_c(\cA) \otimes H^*_c(\cA) \longrightarrow H^*_c(\cA)
\]
makes $H^*_c(\cA)$ a bigraded left module coalgebra over the Hopf algebra $W_0 H^*_c(\cA)$. Both $W_0 H^*_c(\cA)$ and $H^*_c(\cA)$ are connected, and the coaugmentation class of $H^*_c(\cA)$ is $1$ in $H_c^0(\cA_0) \cong \QQ$.
\end{corollary}

Let $\varepsilon_\omega \colon \cH^\omega \to \QQ$ be the augmentation onto the genus-zero summand, and define
\[
j := (\id_{\cH} \otimes \varepsilon_\omega) \rho \colon \cH^\omega \longrightarrow \cH.
\]
The map $j$ is induced on $E^1$ by the forgetful map $f \colon X \to BK(\ZZ)$. Indeed, collapsing the second factor of $\rho_X = (f \times \id_X) \Phi_X$ gives $f \pi_1 \Phi_X$, which is filtered homotopic to $f$. In particular, $j$ is an algebra map.

Since $\bar \rho = (q \otimes \id_{\cH^\omega}) \rho$, applying the augmentation to its second factor gives
\[
(\id_{\overline{\cH}} \otimes \varepsilon_\omega) \bar \rho = qj.
\]
Since the coaugmentation class $1$ is dual to $\varepsilon_\omega$, the map $\iota(a) = \mu(a \otimes 1)$ satisfies
\[
\iota = (qj)^\vee \colon \overline{\cH}^{\vee} \longrightarrow (\cH^\omega)^\vee.
\]

The following consequence of the structure theorem for module coalgebras \cite[Theorem~4.4]{MilnorMoore} reduces freeness to injectivity of $\iota$. The theorem requires an action compatible with the coproduct and counit. It does not require an algebra structure on the module coalgebra or cocommutativity of the acting Hopf algebra. 

For a Hopf algebra $A$ over $\QQ$ and a left $A$-module coalgebra $B$, write $$\cQ_A(B) := \QQ \otimes_A B.$$ This quotient has a natural coalgebra structure, and the projection $\pi \colon B \to \cQ_A(B)$ is a coalgebra map \cite[Proposition~4.3]{MilnorMoore}.

\begin{proposition}\label{prop:MM-freeness}
Let $A$ be a connected graded Hopf algebra over $\QQ$ and $B$ a connected graded left $A$-module coalgebra, with coaugmentation class $1 \in B_0$. If the map ${\iota \colon A \to B}$, $a \mapsto a \cdot 1$, is injective, then there is a graded isomorphism
\[
B \cong A \otimes \cQ_A(B)
\]
of left $A$-modules and right $\cQ_A(B)$-comodules. Thus $B$ is free over $A$. If the structure maps preserve a second grading, the isomorphism can be chosen bigraded.
\end{proposition}

\begin{proof}
Because $1$ is grouplike and the coalgebra structure maps of $B$ are $A$-linear, $\iota$ is a morphism of left $A$-modules and of coalgebras. The projection $\pi \colon B \to \cQ_A(B)$ is surjective, and $\iota$ is injective by hypothesis. Over $\QQ$, the sequences $0 \to A \xrightarrow{\iota} B$ and $B \xrightarrow{\pi} \cQ_A(B) \to 0$ therefore split as sequences of graded vector spaces. Then \cite[Theorem~4.4]{MilnorMoore} gives the asserted graded isomorphism. If the structures are bigraded, the splittings may be chosen one bidegree at a time, and the proposition follows.
\end{proof}

\begin{remark}
The conclusion of Proposition~\ref{prop:MM-freeness} is reminiscent of the structure theorem for Hopf modules, which requires a compatible coaction in addition to the Hopf algebra action. Here freeness instead follows from the compatible coalgebra structure on $B$, connectedness, and injectivity of $\iota$.
\end{remark}

We will also use the following criterion for injectivity of coalgebra maps.

\begin{lemma}\label{lem:coalgebra-injectivity}
A morphism of connected graded coalgebras over a field is injective if and only if it is injective on the subspace of primitive elements.
\end{lemma}

\begin{proof}
Consider a morphism of coalgebras that is injective on primitives. Compatibility with the counit gives injectivity in degree zero. We proceed by induction on degree. Assume injectivity in degrees less than $k$, and let $x$ be a homogeneous element of degree $k$ in the kernel. Write its reduced coproduct as
\[
\overline{\Delta}(x) = \sum_{i = 1}^{k - 1} \sum_j a_{i, j} \otimes b_{k - i, j},
\]
where $a_{i, j}$ has degree $i$ and $b_{k - i, j}$ has degree $k - i$. Applying the coalgebra map to both factors sends this sum to zero. By induction, the induced map on the direct sum of these tensor products is injective, so $\overline{\Delta}(x) = 0$. Thus $x$ is primitive, so $x = 0$.
\end{proof}

\section{Toroidal compactifications and weight-zero cohomology}\label{sec:toroidal-resolution}

In this section, we construct a space $X^\Delta$ with boundary $\partial X^\Delta$ that admits three maps of pairs
\begin{equation}\label{eq:XDelta-maps}
\xymatrix@C=8pt@R=28pt{
& (X^\Delta, \partial X^\Delta) \ar[dl]_{\alpha} \ar[d]^{\beta} \ar[dr]^{\tau} & \\
(A_g^\Delta, \partial A_g^\Delta) & (F_g X, F_{g - 1} X) & (\overline \Sigma^\Delta, \partial \overline \Sigma^\Delta).
}
\end{equation}
Here $A_g^\Delta$ is a toroidal compactification of $A_g$, and $\overline \Sigma^\Delta$ is its extended tropicalization. In this section, we study the maps induced by $\alpha$ and $\tau$ on relative rational cohomology. Pullback along $\alpha$ is an isomorphism, and under the canonical identification $H_c^*(A_g^\trop) \cong W_0 H_c^*(A_g)$, the composite $(\alpha^*)^{-1}\tau^*$ is the natural weight-zero inclusion.

In Section~\ref{sec:injectivity}, we define $\beta$ and compare its pullback with those of $\alpha$ and $\tau$. These comparisons prove Theorems~\ref{thm:free-module} and~\ref{thm:iota-weight-zero}, identifying the map $\iota \colon W_0 H_c^*(\cA) \to H_c^*(\cA)$ induced by the action on the coaugmentation class with the weight-zero inclusion.

\subsection{Toroidal compactifications}\label{sec:toroidal-background}\label{sec:AMRT-tropicalizations}\label{sec:neat-level-covers}

We recall the theory of toroidal compactifications of locally symmetric varieties, due to Ash, Mumford, Rapoport, and Tai. We refer to \cite{AMRT} for details and proofs omitted here. Fix $g$ for the rest of this section. Let $\Gamma := \Sp_{2g}(\ZZ)$, and let $\fH$ denote the Siegel upper half-space. Let $A_g = \Gamma \backslash \fH$ denote the coarse moduli space of $\cA_g = [\Gamma \backslash \fH]$. The natural map $\cA_g \to A_g$ induces an isomorphism on rational cohomology.

Let $P_g^\rt$ be the set of positive semidefinite symmetric $g \times g$ real matrices with rational kernel. A rational polyhedral fan with support $P_g^\rt$ is \emph{admissible} if it is invariant under $\GL_g(\ZZ)$ and has finitely many orbits of cones. It is \emph{smooth} if each cone is generated by a subset of a $\ZZ$-basis of the lattice of symmetric integer matrices. Fix a smooth admissible fan $\Delta$. Give $P_g^\rt$ the colimit topology
\[
P_g^\rt \cong \varinjlim_{\sigma \in \Delta} \sigma,
\]
where each cone has its Euclidean topology and the maps are inclusions of faces. Thus a subset is closed precisely when its intersection with every cone is closed. This topology is independent of the choice of admissible fan, since any two have a common refinement, cf. \cite[Remark~3.2]{BBCMMW}.

Let $I \subseteq \ZZ^{2g}$ be a saturated isotropic sublattice of rank $r$, and let $P_I^\rt \subseteq \Sym^2(I_\RR)$ be the cone of positive semidefinite forms on $I_\RR^\vee$ with rational kernel. A basis of $I$ identifies $P_I^\rt$ with $P_r^\rt$, which we regard as a face of $P_g^\rt$ via $Q \mapsto \operatorname{diag}(Q, 0)$. Restricting $\Delta$ to this face induces a fan $\Delta_I$ with support $P_I^\rt$. It is independent of the choice of basis, since $\Delta$ is $\GL_g(\ZZ)$-invariant. Then $\Delta^\omega := \bigcup_I \Delta_I$, over saturated isotropic sublattices $I$, is a $\Gamma$-invariant fan in $\Sym^2(\RR^{2g})$. 

\begin{definition}\label{def:I-sigma}
For a cone $\sigma \in \Delta^\omega$, let $I_\sigma$ be the smallest saturated isotropic sublattice $I \subseteq \ZZ^{2g}$ such that $\sigma \in \Delta_I$.
\end{definition}

\noindent Note that $I_\tau \subseteq I_\sigma$ for $\tau \preceq \sigma$, and $I_{\gamma \sigma} = \gamma I_\sigma$, for $\gamma \in \Gamma$. Also, $\Delta$ is the restriction of $\Delta^\omega$ to $\Sym^2(\RR^g)$ for the standard isotropic subspace $\RR^g \subseteq \RR^{2g}$, so each of these two fans determines the other.

We write $A_g^\Delta$ for the associated toroidal compactification of $A_g$ and $\partial A_g^\Delta := A_g^\Delta \smallsetminus A_g$ for its boundary. The strata correspond to $\GL_g(\ZZ)$-orbits of cones in $\Delta$ or, equivalently, $\Gamma$-orbits of cones in $\Delta^\omega$.  The space $A_g^\Delta$ is the $\Gamma$-quotient of a bordification $\fH^\Delta$ of $\fH$ whose strata correspond to cones in $\Delta^\omega$. 

From this construction, one can also see that $A_g^\Delta$ is the coarse quotient of a simple normal crossings compactification of a finite cover of $\cA_g$.  Indeed, let $\Gamma' := \Gamma_g(3)$ be the principal congruence subgroup of level 3, and let $G := \Gamma/\Gamma'$. Then $\Delta$ determines a  compactification $A'$ of $\Gamma' \backslash \fH$, and $(A', \partial A')/G = (A_g^\Delta, \partial A_g^\Delta)$. The unipotent translation lattices for $\Gamma'$ are three times those for $\Gamma$, so the fans remain smooth. Moreover, the primitive ray generators of each smooth cone are distinct modulo $3$, whereas $\Gamma'$ acts trivially modulo $3$. Thus no cone contains two distinct $\Gamma'$-equivalent rays, and the boundary of $A'$ is a divisor with simple normal crossings \cite[Proposition~A.3]{MaSiegelOperators}. By construction, the strata of $A'$ correspond to $\Gamma'$-orbits of cones in $\Delta^\omega$.

The tropicalization of $A_g$ with respect to the toroidal compactification $A_g^\Delta$ is 
\[
\Sigma^\Delta = \Delta/\GL_g(\ZZ) = \Delta^\omega/\Gamma.
\]
Its underlying space is canonically identified with $A_g^\trop=P_g^\rt/\GL_g(\ZZ)$, which does not depend on $\Delta$. We will also work with the extended tropicalization, which does depend on $\Delta$, as we now recall.

Let $v_1, \ldots, v_n$ be the primitive generators of the rays of a cone $\sigma \in \Delta$.  We naturally identify $\sigma \cong [0,\infty)^n $ via  $a_1 v_1 + \cdots + a_n v_n \mapsto (a_1, \ldots, a_n)$.  The \emph{extended cone} associated to $\sigma$ is $\overline \sigma \cong [0, \infty]^n$ \cite[\S2.3]{ACP}. For each face $\tau \preceq \sigma$, identify $\overline\tau$ with the closure of $\tau$ in $\overline\sigma$. Gluing these extended cones and taking the quotient gives the extended tropicalization
\[
\overline\Sigma^\Delta
= \Bigl(\varinjlim_{\sigma\in\Delta}\overline\sigma\Bigr)/\GL_g(\ZZ)
= \Bigl(\varinjlim_{\sigma\in\Delta^\omega}\overline\sigma\Bigr)/\Gamma.
\]
It is a compact Hausdorff space that contains $\Sigma^\Delta$ as a dense open subset. Write $$\partial\overline\Sigma^\Delta=\overline\Sigma^\Delta\smallsetminus\Sigma^\Delta.$$  The analogous constructions for the compactification $A'$ of $\Gamma' \backslash \fH$ give the tropicalization $\Sigma(A',\partial A') = \Delta^\omega/\Gamma'$ and extended tropicalization
\[
\overline\Sigma(A',\partial A') = \Bigl(\varinjlim_{\sigma\in\Delta^\omega}\overline\sigma\Bigr)/\Gamma'.
\]
Taking quotients by $G$, we then have
\[
\Sigma(A',\partial A')/G\cong\Sigma^\Delta
\qquad \mbox{ and } \qquad 
\overline\Sigma(A',\partial A')/G\cong\overline\Sigma^\Delta.
\]

\subsection{The Kato--Nakayama spaces}\label{sec:log-cover}

For topological purposes, it is convenient to replace the toroidal pairs $(A_g^\Delta,\partial A_g^\Delta)$ and $(A',\partial A')$ by the associated Kato--Nakayama spaces with their boundaries. The construction removes the simple normal crossing boundary divisor and substitutes the boundary of a manifold with corners. We briefly recall the general construction, which is also known as a real-oriented blowup.

\medskip

Let $Y$ be a smooth complex variety and $D\subseteq Y$ a divisor with simple normal crossings. Near each point of $D$, there are local coordinates $z_1,\ldots,z_n$ on a complex polydisc in which $D$ is given by $z_1\cdots z_r=0$. The real oriented blowup replaces each of the first $r$ coordinate discs by a cylinder
\[
[0,\varepsilon)\times S^1\longrightarrow\{z\in\CC:|z|<\varepsilon\},
\qquad
(r,e^{i\theta})\longmapsto r e^{i\theta}.
\]
In other words, when passing to the boundary of the log space, we remember the angle of approach to each branch of the boundary divisor $D$. These local constructions glue to give the global log space $Y^{\log}$, a manifold with corners equipped with a proper map \cite[\S1.2]{KatoNakayama}
\[
\pi_Y\colon Y^{\log}\longrightarrow Y.
\]
Its boundary is $\pi_Y^{-1}(D)$, and $\pi_Y$ identifies its interior with $Y\smallsetminus D$. Over each locally closed stratum of codimension $r$ in $Y$, the map $\pi_Y$ is an $(S^1)^r$-bundle, with one $S^1$-factor for each branch of $D$ containing the stratum. If $Y$ is proper, then pullback gives the canonical identifications
\[
H^*(Y,D)
\xrightarrow{\;\sim\;}
H^*(Y^{\log},\partial Y^{\log})
\cong H_c^*(Y\smallsetminus D).
\]

Applying this construction to $(A',\partial A')$ gives a compact manifold with corners $(A')^{\log}$, with interior $\Gamma'\backslash\fH$. Write $\pi\colon(A')^{\log}\to A'$ for the projection. The action of $G$ lifts naturally, and $\pi$ is $G$-equivariant. Set
\[
(A_g^\Delta)^{\log}:=(A')^{\log}/G,
\qquad
\partial(A_g^\Delta)^{\log}:=\partial(A')^{\log}/G.
\]
The induced projection to $A_g^\Delta$ identifies the interiors with $A_g$ and induces an isomorphism on relative cohomology.

Since the open inclusion $\Gamma' \backslash \fH \hookrightarrow (A')^{\log}$ is a homotopy equivalence, the universal covering $\fH \to \Gamma' \backslash \fH$ extends to a universal covering
\[
\fH^{\log} \longrightarrow (A')^{\log}
\]
with deck transformation group $\Gamma'$. The space $\fH^{\log}$ is a manifold with corners with interior $\fH$. The action of $\Gamma$ on $\fH$ extends uniquely to $\fH^{\log}$.
Note that the strata of $\fH^{\log}$ correspond naturally to the cones in $\Delta^\omega$. We write $\fH_\sigma^\circ$ for the locally closed stratum of codimension $\dim \sigma$ that corresponds to $\sigma$ and $\fH_\sigma$ for its closure. Then 
\begin{equation} \label{eq:log-strata}
\fH_\sigma = \bigsqcup_{\sigma \preceq \eta} \fH_\eta^\circ.
\end{equation}

\subsection{The space \texorpdfstring{$X^\Delta$}{XDelta}}\label{sec:XgDelta}\label{sec:logarithmic-comparison}

We begin by considering the simplicial space $Y_\bullet^\Delta$ given by 
\[
Y_p^\Delta := \coprod_{\sigma_0 \preceq \cdots \preceq \sigma_p} \fH_{\sigma_p},
\]
with the coproduct topology. The face maps delete cones in the chain $\sigma_0 \preceq \cdots \preceq \sigma_p$ and the degeneracies repeat cones. When we delete the last cone, we also use the inclusion $\fH_{\sigma_p} \subseteq \fH_{\sigma_{p - 1}}$ that follows from \eqref{eq:log-strata}. The chains with $\sigma_0 \ne 0$ form a simplicial subspace $\partial Y_\bullet^\Delta$. Note that $\Gamma$ acts naturally on the pair $(Y_\bullet^\Delta, \partial Y_\bullet^\Delta)$.

\begin{definition}\label{def:XgDelta}
The pair $(X^\Delta, \partial X^\Delta)$ is given by
\[
(X^\Delta, \partial X^\Delta) := E \Gamma \times_\Gamma \bigl(|Y_\bullet^\Delta|, |\partial Y_\bullet^\Delta|\bigr).
\]
\end{definition}

\noindent The $\Gamma$-equivariant forgetful map $|Y_\bullet^\Delta| \to \fH^{\log}$ induces a map $X^\Delta \to (A_g^\Delta)^{\log}$. Composing with the projection from the log space gives the map of pairs $\alpha$ in~\eqref{eq:XDelta-maps}.

Let $\sigma \in \Delta^\omega$ be a cone with primitive ray generators $v_1, \ldots, v_n$. Recall the identification $\overline \sigma \cong [0, \infty]^n$. For each face $\tau \preceq \sigma$, let $\overline v_\tau \in \overline \sigma$ be the point with coordinates
\[
(\overline v_\tau)_i = \begin{cases}
\infty, & v_i \in \tau, \\
0, & v_i \notin \tau.
\end{cases}
\]
Applying $u \mapsto u/(1 + u)$ in each coordinate, with $\infty \mapsto 1$, identifies $\overline \sigma$ with the unit cube $[0, 1]^n$. The standard triangulation of this cube has a simplex with vertices $\overline v_{\sigma_0}, \ldots, \overline v_{\sigma_p}$ for each chain $\sigma_0 \preceq \cdots \preceq \sigma_p \preceq \sigma$ \cite[\S2.5]{ACP}.

We now define the map of pairs
\[
\tau\colon (X^\Delta,\partial X^\Delta)
\longrightarrow
(\overline\Sigma^\Delta,\partial\overline\Sigma^\Delta).
\]
We first give the construction on $|Y_\bullet^\Delta|$. For a chain $\sigma_0\preceq\cdots\preceq\sigma_p$, let $e_0,\ldots,e_p$ denote the vertices of the standard $p$-simplex. For $y \in \fH_{\sigma_p}$, send $(e_i, y)$ to $\overline v_{\sigma_i}$ for $0 \leq i \leq p$, and extend affine linearly. These maps forget $y$, agree on faces, and respect degeneracies. The construction is $\Gamma$-equivariant, so composing with the quotient to $\overline\Sigma^\Delta$ and forgetting the $E\Gamma$ factor gives $\tau$.

For $\Gamma'$, set
\[
(X', \partial X') := E \Gamma \times_{\Gamma'} \bigl(|Y_\bullet^\Delta|, |\partial Y_\bullet^\Delta|\bigr),
\]
and write $\overline\Sigma' := \overline\Sigma(A', \partial A')$. The constructions above give $G$-equivariant maps of pairs
\[
(A', \partial A') \xleftarrow{\alpha'} (X', \partial X') \xrightarrow{\tau'} (\overline\Sigma', \partial\overline\Sigma').
\]
Taking quotients by $G$ gives $(X', \partial X')/G = (X^\Delta, \partial X^\Delta)$ and the maps $\alpha$ and $\tau$.

\begin{proposition}\label{thm:relative-toroidal-resolution}
The maps $\alpha$ and $\alpha'$ induce isomorphisms on relative rational cohomology, and $(\alpha')^*$ is $G$-equivariant.
\end{proposition}

\begin{proof}
Forgetting the chain gives a $\Gamma$-equivariant map of pairs
\[
e_{\log} \colon (|Y_\bullet^\Delta|, |\partial Y_\bullet^\Delta|)
\longrightarrow (\fH^{\log}, \partial \fH^{\log}).
\]
Choose a triangulation of $(A')^{\log}$ in which the closed strata are subcomplexes, and lift it to $\fH^{\log}$. Over an open simplex contained in $\fH_\eta^\circ$, evaluation is the projection from the product with the realization of the poset of faces of $\eta$. For the boundary, use the nonzero faces. Both posets have terminal element $\eta$, so their realizations are contractible. Filtering by the skeleta of the triangulation therefore shows that $e_{\log}$ induces homology isomorphisms on the spaces and their boundaries.

The Serre spectral sequence over $B\Gamma'$ then gives an isomorphism on relative rational cohomology after taking homotopy quotients. Since $\Gamma'$ acts freely and properly, projection from the homotopy quotient to $\Gamma' \backslash \fH^{\log} = (A')^{\log}$ is a homotopy equivalence of pairs. By Section~\ref{sec:log-cover}, the projection $(A')^{\log} \to A'$ induces an isomorphism on relative cohomology. Thus $(\alpha')^*$ is a $G$-equivariant isomorphism. Since $G$ is finite, taking $G$-invariants proves the assertion for $\alpha$.
\end{proof}

\subsection{Tropical pullback and weight-zero cohomology}\label{sec:incidence-weight-zero}

Recall the canonical isomorphism
\[
\vartheta \colon H_c^*(A_g^\trop) \xrightarrow{\ \sim \ } W_0 H_c^*(A_g)
\]
of \cite[Proposition~2.1]{BCGP}. We will show that the maps $\alpha$ and $\tau$ realize this comparison, after identifying compactly supported cohomology of $A_g$ and $A_g^\trop$ with the cohomology of the pairs $(A_g^\Delta, \partial A_g^\Delta)$ and $(\overline\Sigma^\Delta, \partial\overline\Sigma^\Delta)$, respectively.

\begin{proposition}\label{prop:incidence-weight-zero}
The composite
\[
(\alpha^*)^{-1}\tau^* \colon H_c^*(A_g^\trop) \longrightarrow H_c^*(A_g)
\]
is $\vartheta$ followed by the natural inclusion $W_0 H_c^*(A_g) \hookrightarrow H_c^*(A_g)$.
\end{proposition}

\begin{proof}
We work $G$-equivariantly with the pair $(A',\partial A')$, following the weight spectral sequence argument in \cite[Proposition~4.3]{HainautPetersen}. We will then deduce the result for $A_g$ by passing to $G$-invariants.

The simplicial resolution of $(A',\partial A')$ indexed by chains of connected closed strata gives a $G$-equivariant spectral sequence of mixed Hodge structures \cite[\S3]{PetersenStratified}
\[
E_1^{p,q}=\bigoplus_{A'=Z_0\supsetneq\cdots\supsetneq Z_p}H^q(Z_p)\Longrightarrow H^{p+q}(A',\partial A').
\]
The closed strata are smooth and proper, so row $q$ is pure of weight $q$. Consequently, the spectral sequence degenerates at $E_2$, and
\[
E_2^{p,0}\cong W_0H^p(A',\partial A').
\]
The row $q=0$ is the relative simplicial cochain complex of $(\overline\Sigma',\partial\overline\Sigma')$. Taking cohomology gives a $G$-equivariant isomorphism
\[
\vartheta'\colon H^*(\overline\Sigma',\partial\overline\Sigma')\xrightarrow{\sim}W_0H^*(A',\partial A').
\]

The projection $\fH^{\log} \to A'$ induces a map from $X'$ to the realization of the simplicial resolution. Composing with evaluation to $A'$ gives $\alpha'$, while forgetting the point in each stratum gives $\tau'$. Pullback along the projection to $\overline\Sigma'$ identifies its relative simplicial cochain complex with the $q=0$ row, whose inclusion induces $\vartheta'$. It follows that $(\tau')^* = (\alpha')^* \vartheta'$. By naturality, the map induced by $\vartheta'$ on $G$-invariants is $\vartheta$. Taking $G$-invariants therefore gives $\tau^* = \alpha^* \vartheta$, and the proposition follows.
\end{proof}

\section{Freeness and the weight-zero inclusion}\label{sec:injectivity}

Recall from Section~\ref{sec:comod} that we have given $H_c^*(\cA)$ the structure of a bigraded module coalgebra over the Hopf algebra $W_0 H_c^*(\cA)$. To prove Theorem~\ref{thm:free-module}, it remains to show that the map
\[
\iota \colon W_0 H_c^*(\cA) \longrightarrow H_c^*(\cA), \qquad a \longmapsto a \cdot 1,
\]
obtained by acting on the coaugmentation class, is injective. We will then prove Theorem~\ref{thm:iota-weight-zero}, which identifies $\iota$ with the natural inclusion of the weight-zero subspace.

\subsection{The map to the symplectic \texorpdfstring{$S$}{S}-construction}\label{sec:flag-map}

We continue to use the notation of Section~\ref{sec:toroidal-resolution}. In particular, $\Delta$ is a smooth admissible fan with support $P_g^\rt$, and $\Delta^\omega = \bigcup_I \Delta_I$ is the associated fan in $\Sym^2(\RR^{2g})$. Then $Y_\bullet^\Delta$ is the simplicial space
\[
Y_p^\Delta = \coprod_{\sigma_0 \preceq \cdots \preceq \sigma_p} \fH_{\sigma_p},
\]
where the coproduct ranges over chains of cones in $\Delta^\omega$, and $X^\Delta = E \Gamma \times_\Gamma |Y_\bullet^\Delta|$.

For each cone $\sigma \in \Delta^\omega$, let $I_\sigma \subseteq \ZZ^{2g}$ be the smallest saturated isotropic sublattice such that the span of $\Sym^2(I_\sigma)$ contains $\sigma$ (Definition~\ref{def:I-sigma}).  We now construct a map from $X^\Delta$ to the symplectic $S$-construction $X$ from Section~\ref{sec:symplectic-S}. Let $\sigma_\bullet = (\sigma_0 \preceq \cdots \preceq \sigma_p)$ be a chain of cones in $\Delta^\omega$. Set $I_j := I_{\sigma_j}$, so
\[
I_0 \subseteq \cdots \subseteq I_p.
\]
For $0 \leq i \leq p$ and $i < j \leq 2p + 1 - i$, let
\[
Q_{i, j} :=
\begin{cases}
I_j/I_i, & j \leq p,\\[2pt]
I_{2p + 1 - j}^\perp/I_i, & j > p,
\end{cases}
\]
with pairing induced by $\omega$. The natural inclusions and quotient maps give a diagram of the form~\eqref{eq:Qdiagram}, and hence an object $Q(\sigma_\bullet)$ of $\cC_p$.

\begin{proposition}\label{prop:flag-map}
Sending $(\sigma_\bullet, y)$ to $Q(\sigma_\bullet)$, and each element $\gamma \in \Gamma$ to the induced isomorphism of diagrams, defines a map of pairs
\[
\beta \colon \bigl(X^\Delta, \partial X^\Delta\bigr) \longrightarrow \bigl(F_g X, F_{g - 1} X\bigr).
\]
\end{proposition}

\begin{proof}
The assignment $\sigma_\bullet \mapsto Q(\sigma_\bullet)$ is $\Gamma$-equivariant and gives continuous maps
\[
\Gamma^q \times Y_p^\Delta \longrightarrow N_q \cC_p
\]
that assemble to a map of bisimplicial spaces. Taking realizations gives the desired map $\beta$. The symplectic lattice $I_0^\perp/I_0$ has genus $g - \rank I_0$, so the image lies in $F_g X$. On $\partial X^\Delta$, we have $\sigma_0 \neq 0$ and hence $I_0 \neq 0$, so the image of $\partial X^\Delta$ lies in $F_{g - 1} X$.
\end{proof}

\subsection{Comparison with tropicalization}\label{sec:tropical-comparison}

Recall that $f \colon X \to BK(\ZZ)$ is the forgetful map that sends an isotropic flag in a symplectic lattice to its underlying flag of lattices, as defined in Section~\ref{sec:comod}. In this section, we compare the map on relative cohomology induced by $f \circ \beta$ with the map induced by tropicalization $\tau \colon X^\Delta \to \overline \Sigma^\Delta$.

Let $U_g := A_g^\trop \smallsetminus A_{g - 1}^\trop$ be the open subset of forms of rank $g$, with $A_{-1}^\trop := \emptyset$. Write $U_g^+$ for its one-point compactification, with distinguished point $*$. We identify $H_c^*(U_g)$ with $H^*(U_g^+, *)$. The construction in \cite[\S2.4]{BCGP} gives an isomorphism
\[
\vartheta_Q \colon H^*(U_g^+, *) \xrightarrow{\ \sim \ } H^*(F_g BK(\ZZ), F_{g - 1} BK(\ZZ))
\]
whose construction we recall in the proof. Since $U_g$ is an open subset of $\overline \Sigma^\Delta \smallsetminus \partial \overline \Sigma^\Delta$, collapsing its complement gives a map of pairs
\[
c \colon (\overline \Sigma^\Delta, \partial \overline \Sigma^\Delta) \longrightarrow (U_g^+, *).
\]
Under the identification above, $c^*$ is extension by zero.

\begin{proposition}\label{prop:bcgp-incidence-square}
The following diagram commutes:
\[
\begin{tikzcd}[column sep=3em,row sep=1.5em]
& H^*(F_g BK(\ZZ), F_{g - 1} BK(\ZZ)) \arrow[d,"f^*"] \\
H^*(U_g^+, *) \arrow[ur,"\vartheta_Q"] \arrow[d,"c^*"'] & H^*(F_g X, F_{g - 1} X) \arrow[d,"\beta^*"] \\
H^*(\overline \Sigma^\Delta, \partial \overline \Sigma^\Delta) \arrow[r,"\tau^*"'] & H^*(X^\Delta, \partial X^\Delta) \ .
\end{tikzcd}
\]
\end{proposition}

In preparation for the proof of Proposition~\ref{prop:bcgp-incidence-square}, we choose, for each cone $\sigma \in \Delta^\omega$, a finite set of primitive vectors $V_\sigma \subseteq I_\sigma$, closed under negation and spanning $I_\sigma \otimes \QQ$, compatibly with face inclusions and the action of $\Gamma$.

\begin{lemma}\label{lem:compatible-spanning-sets}
There is a $\Gamma$-equivariant family $\{ V_\sigma : \sigma \in \Delta^\omega \}$ of finite sets of primitive vectors $V_\sigma \subseteq I_\sigma$ spanning $I_\sigma \otimes \QQ$, such that $V_\sigma = -V_\sigma$ and $V_\tau \subseteq V_\sigma$ whenever $\tau \preceq \sigma$.
\end{lemma}

\begin{proof}
The stabilizer in $\Gamma$ of a ray $\rho$ fixes its primitive integral generator, a positive definite form on $I_\rho^\vee \otimes \RR$, so its image in $\GL(I_\rho)$ is finite. The images of basis vectors and their negatives under this finite group give an invariant spanning set $V_\rho$ closed under negation. Choose such sets for one ray in each $\Gamma$-orbit and extend equivariantly. Set $V_\sigma := \bigcup_{\rho \preceq \sigma, \, \dim \rho = 1} V_\rho$. Equivariance, closure under negation, and face compatibility follow from the construction. Since the inclusions $I_\rho \subseteq I_\sigma$ are saturated and the rays generate $\sigma$, these sets consist of primitive vectors and span $I_\sigma \otimes \QQ$.
\end{proof}

Fix such a family $\{ V_\sigma : \sigma \in \Delta^\omega \}$.

\begin{proof}[Proof of Proposition~\ref{prop:bcgp-incidence-square}]
For $g = 0$, all maps in the diagram are identities. Assume $g > 0$. We recall the construction of the isomorphism $\vartheta_Q$ from \cite[\S2.4]{BCGP}. Let $Z := |N_\bullet T_\bullet(\QQ^g) \cup \{*\} |$, with distinguished point $*$. The maps in that construction fit into the following diagram of pairs:
\begin{equation}\label{eq:beta}
\begin{tikzcd}[column sep=3em,row sep=1.5em]
(Z, *) \arrow[r,"a"] \arrow[d,"b"'] & (F_g BK(\ZZ)/F_{g - 1} BK(\ZZ), *) \\
(U_g^+, *) & (F_g X/F_{g - 1} X, *) \arrow[u,"f"'] \\
(\overline \Sigma^\Delta, \partial \overline \Sigma^\Delta) \arrow[u,"c"] & (X^\Delta, \partial X^\Delta) \arrow[u,"\beta"'] \arrow[l,"\tau"]
\end{tikzcd}
\end{equation}
Here $a$ is the map of \cite[\S2.4]{BCGP} followed by duality and reflection. The maps $a$ and $b$ induce isomorphisms on reduced rational homology. The isomorphism $\vartheta_Q$ is $(a^*)^{-1} b^*$. Thus the proposition asserts that the two routes from $H^*(U_g^+, *)$ to $H^*(X^\Delta, \partial X^\Delta)$ in the cohomology of \eqref{eq:beta} agree, using $(a^*)^{-1}$ to traverse the zig-zag.

We prove this equality by constructing a chain map $\ell$ that lifts relative rational chains of $(X^\Delta, \partial X^\Delta)$ to those of $(Z, *)$ such that
\begin{equation} \label{eq:lift}
a_* \ell = (f \circ \beta)_*, \qquad \mbox{ and } \qquad b_* \ell \simeq (c \circ \tau)_*,
\end{equation}
where $\simeq$ denotes a chain homotopy. The proposition follows from the existence of such a lift, by dualizing and passing to cohomology.

Let $(\sigma_\bullet, y) \in Y_p^\Delta$ be a $p$-simplex, with $\sigma_\bullet = (\sigma_0 \preceq \cdots \preceq \sigma_p)$. If $\sigma_0 \neq 0$ or $\rank I_p < g$, set $\ell([\sigma_\bullet, y]) = 0$. Otherwise, let $V_i := V_{\sigma_i}$ be the fixed finite set of primitive vectors in $I_i$,  chosen as in Lemma~\ref{lem:compatible-spanning-sets}, and let $N := |V_p |$. Each total order $<$ on $V_p$, restricted to the subsets $V_i$, determines an oriented $p$-simplex $[V_p \supseteq \cdots \supseteq V_0, <]$ in $Z$. Define $\ell$ by
\begin{equation} \label{eq:ell}
\ell([\sigma_\bullet, y]) := \frac{(-1)^{p(p + 1)/2}}{N!} \sum_<[V_p \supseteq \cdots \supseteq V_0, <],
\end{equation}
where the sum ranges over all total orders on $V_p$. The sign accounts for orientations. Each order on a subset of $m$ elements has $N!/m!$ extensions, so averaging commutes with restriction. The face and degeneracy maps agree with those in the construction of $Z$, with the order of the faces reversed. The displayed sign therefore makes $\ell$ commute with the differential. The assignment is natural under lattice isomorphisms, so it also commutes with the nerve operators and defines the required map on chains of $X^\Delta$.

It remains to show that $\ell$ satisfies the two identities in \eqref{eq:lift}. For the first identity, write $M := I_p^\vee$ and let $V_i^\perp \subseteq M$ be the common kernel of the covectors in $V_i$. The diagram associated to the ordered sets $V_p \supseteq \cdots \supseteq V_0$ has modules
\[
P_{ij} = V_{p - j}^\perp/V_{p - i}^\perp, \qquad (P_{p - j, p - i})^\vee \cong I_j/I_i,
\]
where the second identification follows because $V_i$ spans $I_i \otimes \QQ$ and the $I_j$ are saturated. Thus dualizing and reflecting this diagram gives the diagram defining $f \circ \beta$. The reflection also reverses the simplex coordinates, canceling the sign in \eqref{eq:ell}. It follows that $a_* \ell = (f \circ \beta)_*$.

For the second identity, we construct a map of pairs
\[
h \colon (X^\Delta, \partial X^\Delta) \longrightarrow (U_g^+, *)
\]
inducing $b_* \ell$, and a homotopy of pairs from $c \circ \tau$ to $h$.

Consider a simplex indexed by $\sigma_\bullet$ with $\sigma_0 = 0$ and $\rank I_p = g$. Let $s_0, \ldots, s_p$ be its barycentric coordinates and set
\[
c_j := \frac{s_j + \cdots + s_p}{s_0 + \cdots + s_{j - 1}}, \qquad 1 \leq j \leq p.
\]
For $1 \leq j \leq p$, define positive semidefinite quadratic forms $\varphi_j$ and $\varphi'_j$ on $M_\RR$ as follows. Let $\varphi_j$ be the sum of the primitive generators of the rays of $\sigma_j$ that are not rays of $\sigma_{j - 1}$, and set
\[
\varphi'_j := \sum_{v \in V_j \smallsetminus V_{j - 1}} v^2.
\]
On the interiors of these simplices, $c \circ \tau$ is represented by $\sum_j c_j \varphi_j$, while $h$ is represented by $\sum_j c_j \varphi'_j$. With the coordinates reversed as in \eqref{eq:ell}, the latter expression is the formula for $b$ on each lifted simplex, by \cite[formula~(25)]{BCGP} with the function $z \mapsto z/(1 - z)$. It is independent of the total order, so averaging gives the map on chains induced by $h$. We join these maps by the straight-line homotopy
\[
H_u = \sum_{j = 1}^p c_j\bigl((1 - u) \varphi_j + u \varphi'_j\bigr), \qquad 0 \leq u \leq 1.
\]
Thus $H_0 = c \circ \tau$ and $H_1 = h$. If $\sigma_0 \neq 0$ or $\rank I_p < g$, set $H_u = *$ for all $u$.

The partial sums $\sum_{j \leq k} \varphi_j$ and $\sum_{j \leq k} \varphi'_j$ both have kernel $V_k^\perp \otimes \RR$. Hence $H_u$ is positive definite on the interiors of these simplices. On the face $s_k = 0$, for $1 \leq k < p$, we have $c_k = c_{k + 1}$. Deleting $\sigma_k$ replaces $\varphi_k, \varphi_{k + 1}$ by their sum, and similarly for $\varphi'_k, \varphi'_{k + 1}$. Thus the homotopy on each simplex restricts to the homotopy assigned to each of these faces. The same holds on the face $s_p = 0$ when $\rank I_{p - 1} = g$, since $c_p = 0$. It remains to check continuity at the faces mapped to the basepoint, jointly in the simplex coordinates and $u$. Once this continuity is established, $H_u$ is a homotopy of pairs from $c \circ \tau$ to $h$, and hence induces the required chain homotopy.

To check continuity at the faces mapped to the basepoint, it suffices to show that $H_u$ leaves every compact subset of $U_g$, uniformly in $u$. On positive definite forms on $M_\RR$, the two functions
\[
\nu(Q) = \min_{0 \neq x \in M} Q(x), \qquad \nu^\vee(Q) = \min_{0 \neq v \in M^\vee} Q^{-1}(v)
\]
are continuous, positive, and invariant under lattice isomorphisms. They are therefore bounded below on every compact subset of $U_g$.

As we approach a face whose first remaining cone is nonzero, some $c_k$ tends to infinity with $I_k \neq 0$. Fix $0 \neq v \in I_k$. Since the two partial sums have kernel $V_k^\perp \otimes \RR$, their convex combinations dominate $\kappa v^2$ for a fixed $\kappa > 0$, independent of $u$. Since the coefficients $c_j$ are nonincreasing in $j$, we have $\kappa c_k v^2 \leq H_u$, and hence
\[
\nu^\vee(H_u) \leq \frac{1}{\kappa c_k}.
\]
Thus $\nu^\vee(H_u)$ tends to zero uniformly in $u$.

If the last remaining cone of a face is $\sigma_k$ with $\rank I_k < g$, choose $0 \neq x \in V_k^\perp$. For $j \leq k$, the forms $\varphi_j$ and $\varphi'_j$ vanish on $x$. As we approach this face, the remaining coefficients tend to zero, so $\nu(H_u) \leq H_u(x) \to 0$, uniformly in $u$. This shows that $H_u$ extends to a continuous homotopy of pairs from $c \circ \tau$ to $h$. The construction is invariant under lattice isomorphisms, so the homotopy descends to $X^\Delta$, and the proposition follows.
\end{proof}

\subsection{Completion of the proof of Theorem~\ref{thm:free-module}}\label{sec:proof-of-freeness}

In Section~\ref{sec:comod}, we constructed an action $\mu$ making $H_c^*(\cA)$ a connected bigraded module coalgebra over $W_0 H_c^*(\cA)$. By Proposition~\ref{prop:MM-freeness}, it remains to prove injectivity of the action on the coaugmentation class
\[
\iota \colon W_0 H_c^*(\cA) \longrightarrow H_c^*(\cA), \qquad a \longmapsto a \cdot 1.
\]

Fix $g$. By Proposition~\ref{prop:E1X}, followed by Borel--Serre duality and Poincar\'e duality, we have isomorphisms
\begin{equation}\tag{$\dagger$} \label{eq:symplectic-relative-cohomology}
H^m(F_g X, F_{g - 1} X) \cong H_{m - g}\bigl(\Sp_{2g}(\ZZ); \St_g^\omega\bigr)^\vee \cong H_c^m(A_g).
\end{equation}
Recall from Section~\ref{sec:comod} that $q \colon \cH \to \overline \cH = \cH/(e)$ is the quotient map and $$j = (\id \otimes \varepsilon_\omega) \rho \colon \cH^\omega \to \cH$$ is induced by $f \colon X \to BK(\ZZ)$, where $\varepsilon_\omega$ is the augmentation onto the genus-zero summand. The action is dual to $(q \otimes \id) \rho$, and the coaugmentation class $1$ is dual to $\varepsilon_\omega$.  Thus $\iota = (qj)^\vee$. Taking genus $g$ components and using the identification $H_c^*(A_g^\trop) \cong \overline \cH_g^\vee$ of \cite{BCGP}, we see that $\iota_g \vartheta$ is
\[
H_c^*(A_g^\trop) \xrightarrow{\ q_g^\vee \ } H^*(F_g BK(\ZZ), F_{g - 1} BK(\ZZ)) \xrightarrow{\ f^* \ } H^*(F_g X, F_{g - 1} X),
\]
followed by \eqref{eq:symplectic-relative-cohomology}.

Under the identifications of \cite[Proposition~4.4]{BCGP}, the natural map
\[
H_*^{\mathrm{BM}}(A_g^\trop) \longrightarrow H_*^{\mathrm{BM}}(A_g^\trop, A_{g - 1}^\trop)
\]
provides a right inverse to $q_g$. Its dual is extension by zero from $U_g$, represented by $c^* \vartheta_Q^{-1}$, so
\[
c^* \vartheta_Q^{-1} q_g^\vee = \id.
\]
The comparison $\beta^* f^* \vartheta_Q = \tau^* c^*$ of Proposition~\ref{prop:bcgp-incidence-square} now gives
\begin{equation}\label{eq:injective-composite}
\beta^* f^* q_g^\vee = \tau^*.
\end{equation}

Since $\tau^*$ is injective by Proposition~\ref{prop:incidence-weight-zero}, $f^* q_g^\vee$ is injective. Its composite with \eqref{eq:symplectic-relative-cohomology} is $\iota_g \vartheta$, and $\vartheta$ is an isomorphism, so $\iota_g$ is injective, as required. \qed

\subsection{Identification with the weight-zero inclusion}\label{sec:weight-zero-inclusion}

It remains to identify $\iota$ with the natural inclusion of the weight-zero subspace, as asserted in Theorem~\ref{thm:iota-weight-zero}. The case $g = 0$ is trivial, so throughout this subsection we assume $g > 0$.

\medskip

Let $\Gamma = \Sp_{2g}(\ZZ)$ and let $T = T^\omega_g(\QQ)$ be the symplectic Tits building.  We work with the cone over $T$, whose vertex $0 \in \Cone(T)$ corresponds to the isotropic subspace $0 \subseteq \QQ^{2g}$. We will study and compare maps to $\Cone(T)$ from two different bordifications of $\fH$, the Borel--Serre bordification $\fH^{\BS}$ and the logarithmic bordification $\fH^{\log}$ constructed in Section~\ref{sec:log-cover}. Both are manifolds with corners with interior $\fH$. Each boundary stratum is associated to a rational isotropic subspace as follows. A stratum of $\fH^{\BS}$ corresponds to an isotropic flag $0 \subsetneq I_1 \subsetneq \cdots \subsetneq I_k$ and is associated to its largest subspace $I_k$. A stratum of $\fH^{\log}$ corresponds to a cone $\sigma \in \Delta^\omega$ and is associated to $I_\sigma \otimes \QQ$, where $I_\sigma$ is the lattice of Definition~\ref{def:I-sigma}, i.e., $I_\sigma$ is the smallest saturated isotropic sublattice $I$ such that $\sigma \in \Delta_I$.

For each nonzero isotropic subspace $I$, let $T_I \subseteq T$ be the subcomplex of flags of isotropic subspaces that contain $I$. It is a cone with vertex $I$.

\begin{proposition}\label{prop:boundary-maps}
There are $\Gamma$-equivariant maps of pairs $\varphi_{\log}$ and $\varphi_{\BS}$ to $(\Cone(T),T)$ from $(\fH^{\log}, \partial \fH^{\log})$ and $(\fH^{\BS}, \partial \fH^{\BS})$, respectively, with the property that the image of a boundary stratum associated to $I$ lies in $T_I$. Moreover, each of these maps is unique up to a homotopy through maps of pairs with the same property.
\end{proposition}

\begin{proof}
We use the same construction for both bordifications. Choose a locally finite $\Gamma$-equivariant triangulation of the bordification $M$ in which every closed stratum is a subcomplex. For a boundary simplex $\tau \subseteq \partial M$, let $I_\tau$ be the isotropic subspace associated to the stratum containing its relative interior. If $\tau'$ is a face of $\tau$, then $I_\tau \subseteq I_{\tau'}$. Sending the barycenter of $\tau$ to $I_\tau$ therefore defines an equivariant simplicial map from the barycentric subdivision of the boundary to $T$. This extends to an equivariant simplicial map of pairs from $(M, \partial M)$ to $(\Cone(T), T)$ by sending the barycenter of each simplex not contained in $\partial M$ to the cone point $0$.

Any two equivariant maps $\partial M \to T$ with the property that the image of a boundary stratum associated to $I$ lies in $T_I$ are equivariantly homotopic through such maps. Working in the barycentric subdivision, extend the homotopy one orbit of simplices at a time. On a simplex $\tau$, the extension must take values in $T_{I_\tau}$. The homotopies already defined on its faces take values in this cone, since $T_{I_{\tau'}} \subseteq T_{I_\tau}$. Moreover, the stabilizer of $\tau$ fixes $\tau$ pointwise and fixes $I_\tau$, so its fixed point space in $T_{I_\tau}$ contracts to $I_\tau$. Thus the extension exists by \cite[Chapter~II, \S1]{Bredon}.

To prove uniqueness on $M$, let $\varphi_0$ and $\varphi_1$ be two equivariant maps of pairs mapping each boundary stratum associated to $I$ into $T_I$. The preceding boundary homotopy, together with $\varphi_0$ and $\varphi_1$, defines an equivariant map
\[
(M \times \{0, 1\}) \cup (\partial M \times [0, 1]) \longrightarrow \Cone(T).
\]
Extend this map over the prisms $\tau \times [0, 1]$, one orbit of simplices at a time. The stabilizer of $\tau$ fixes the prism pointwise, and its fixed point space in $\Cone(T)$ contracts to $0$, so the extension exists, and the proposition follows.
\end{proof}

We regard the induced pullbacks as maps $\varphi_{\log}^*, \varphi_{\BS}^* \colon H_\Gamma^*(\Cone(T), T) \to H_c^*(A_g)$, using the canonical identification of the equivariant relative cohomology of each bordification with $H_c^*(A_g)$. Write $\kappa = \kappa_g$ for the map of Lemma~\ref{lem:building-filtration}.

\begin{proposition}\label{prop:BS-duality}
The composite $\varphi_{\BS}^* \circ \kappa^*$ agrees with the isomorphism~\eqref{eq:symplectic-relative-cohomology} given by Borel--Serre and Poincar\'e duality.
\end{proposition}

\begin{proposition}\label{prop:log-tropical-pullback}
We have $\alpha^* \circ \varphi_{\log}^* \circ \kappa^* = \beta^*$.
\end{proposition}

\begin{proposition}\label{prop:equal-boundary-pullbacks}
The pullbacks induced by $\varphi_{\log}$ and $\varphi_{\BS}$ agree: $\varphi_{\log}^* = \varphi_{\BS}^*$.
\end{proposition}

Given these three propositions, Theorem~\ref{thm:iota-weight-zero} follows easily. 

\begin{proof}[Proof of Theorem~\ref{thm:iota-weight-zero}]
The description of $\iota_g$ in Section~\ref{sec:proof-of-freeness} and Proposition~\ref{prop:BS-duality} identify $\iota_g \vartheta$ with $\varphi_{\BS}^* \circ \kappa^* \circ f^* \circ q_g^\vee$. Using Propositions~\ref{prop:equal-boundary-pullbacks} and~\ref{prop:log-tropical-pullback}, we obtain
\[
\alpha^* \iota_g \vartheta = \beta^* f^* q_g^\vee = \tau^*,
\]
where the last equality follows from Proposition~\ref{prop:bcgp-incidence-square}, as in~\eqref{eq:injective-composite}. By Proposition~\ref{prop:incidence-weight-zero}, $(\alpha^*)^{-1} \tau^*$ is $\vartheta$ followed by the weight-zero inclusion. Hence $\iota_g$ is that inclusion. Finally, associativity gives $a \cdot \iota(b) = \iota(ab)$ for $a, b \in W_0 H_c^*(\cA)$, so the action extends the Hopf algebra product.
\end{proof}

We now proceed with the proofs of the propositions.

\subsubsection{Proof of Proposition~\ref{prop:BS-duality}}\label{sec:BS-duality}
For a nonzero rational isotropic subspace $I \subseteq \QQ^{2g}$, let $F_I$ be the closed boundary face corresponding to the maximal parabolic stabilizing $I$, and write $F_0 = \fH^{\BS}$. Consider the simplicial space
\[
Y_p^{\BS} = \coprod_{I_0 \subseteq \cdots \subseteq I_p} F_{I_0} \cap \cdots \cap F_{I_p},
\]
where the indices are rational isotropic subspaces, including $0$. Face maps omit a subspace and use the inclusions of intersections; degeneracies repeat a subspace. The summands with $I_0 \ne 0$ form $\partial Y_\bullet^{\BS}$. Evaluation gives a map of pairs
\[
e_{\BS} \colon (|Y_\bullet^{\BS}|, |\partial Y_\bullet^{\BS}|) \longrightarrow (\fH^{\BS}, \partial \fH^{\BS}).
\]
Sending the vertex indexed by $I$ to $I \in \Cone(T)$ and extending linearly on each simplex defines a map of pairs
\[
\beta_{\BS} \colon (|Y_\bullet^{\BS}|, |\partial Y_\bullet^{\BS}|) \longrightarrow (\Cone(T), T).
\]

We first compare $\beta_{\BS}$ with $\varphi_{\BS} \circ e_{\BS}$. Let $x \in |\partial Y_\bullet^{\BS}|$, and suppose $y = e_{\BS}(x)$ lies in the stratum indexed by $J_1 \subsetneq \cdots \subsetneq J_k$. If $I$ and $I'$ are vertices of the simplices containing $\beta_{\BS}(x)$ and $\varphi_{\BS}(y)$ in their relative interiors, respectively, then $I$ belongs to this flag and $J_k \subseteq I'$. Thus $I \subseteq J_k \subseteq I'$, so the vertices of the two simplices together form a flag. Both image points lie in the corresponding simplex of $T$, and the straight line segments joining them give an equivariant homotopy on $|\partial Y_\bullet^{\BS}|$, compatible with faces and degeneracies. Using the triangulation from Proposition~\ref{prop:boundary-maps}, this extends equivariantly over $|Y_\bullet^{\BS}|$, since the fixed point space of every stabilizer in $\Cone(T)$ contracts to $0$. Thus $\beta_{\BS} \simeq_\Gamma \varphi_{\BS} \circ e_{\BS}$ as maps of pairs.

Both $e_{\BS}$ and $\beta_{\BS}$ induce homology isomorphisms on the spaces and on their boundaries. We first consider $e_{\BS}$. Since $F_0 = \fH^{\BS}$, we may adjoin $0$ to each indexing flag. The straight line homotopy to the vertex indexed by $0$, keeping the point of $\fH^{\BS}$ fixed, gives a deformation retraction of $|Y_\bullet^{\BS}|$ onto $F_0$. For the restriction to the boundary, choose a triangulation of $\partial \fH^{\BS}$ compatible with the strata. Over each open simplex, $e_{\BS}$ is the projection from its product with the simplex of closed boundary faces containing it. The filtration by skeleta therefore shows that the boundary restriction also induces an isomorphism on homology. For $\beta_{\BS}$, every intersection indexed by a flag is contractible, and the nerve of the closed boundary faces is $T$ \cite[\S\S8.3--8.4]{BorelSerre}. The realization spectral sequence therefore shows that $\beta_{\BS}$ induces homology isomorphisms on the spaces and on their boundaries. These isomorphisms are preserved on passing to homotopy quotients, and the long exact sequences of the pairs give the corresponding isomorphisms on relative homology.

The relative simplicial chain complex of $(\Cone(T), T)$ is generated by strict flags beginning with $0$. The face obtained by deleting $0$ lies in $T$, so it is zero in the relative complex. The remaining face maps identify this complex with $\Sigma \widetilde C_*(T; \QQ)$, with differential $-d$. Its homology is $\St_g^\omega$ in degree $g$ and vanishes in all other degrees. Taking derived coinvariants with the group resolution first, as in Lemma~\ref{lem:building-filtration}, and using the comparison above gives
\[
H_m^\Gamma(\fH^{\BS}, \partial \fH^{\BS}; \QQ) \xrightarrow{\sim} H_{m - g}(\Gamma; \St_g^\omega).
\]

Write $d = g(g + 1)$ for the real dimension of $A_g$. Composing the rational orbifold Poincar\'e--Lefschetz isomorphism $H^{d - m}(A_g) \cong H_m^\Gamma(\fH^{\BS}, \partial \fH^{\BS}; \QQ)$ with the displayed map gives Borel--Serre duality \cite[\S8.6 and Theorem~11.4.2]{BorelSerre}. Thus the dual of the displayed map is the second isomorphism in~\eqref{eq:symplectic-relative-cohomology}. Since $\beta_{\BS} \simeq_\Gamma \varphi_{\BS} \circ e_{\BS}$, this dual map is induced by $\varphi_{\BS}$. Finally, Lemma~\ref{lem:building-filtration} shows that $\kappa^*$ gives the first isomorphism in~\eqref{eq:symplectic-relative-cohomology} under the identification $H_\Gamma^m(\Cone(T), T) \cong H_{m - g}(\Gamma; \St_g^\omega)^\vee$ obtained above. Hence $\varphi_{\BS}^* \circ \kappa^*$ is the isomorphism~\eqref{eq:symplectic-relative-cohomology}. \qed

\subsubsection{Proof of Proposition~\ref{prop:log-tropical-pullback}}\label{sec:incidence} Let $e_{\log}$ be the evaluation map from Section~\ref{sec:XgDelta}. Define a flag map $\beta_{\log}$ from $(|Y_\bullet^\Delta|, |\partial Y_\bullet^\Delta|)$ to $(\Cone(T), T)$ by sending the vertex indexed by $\sigma_j$ to $I_{\sigma_j} \otimes \QQ$ and extending linearly on each simplex. The zero cone maps to $0$, and a flag of nonzero cones maps into $T$.

Let $x \in |\partial Y_\bullet^\Delta|$, with $e_{\log}(x) = y \in \fH_\eta^\circ$. If $I$ and $I'$ are vertices of the simplices containing $\beta_{\log}(x)$ and $\varphi_{\log}(y)$ in their relative interiors, respectively, then
\[
I \subseteq I_\eta \otimes \QQ \subseteq I'.
\]
Thus the vertices of these two simplices together form a flag. Both image points lie in the corresponding simplex of $T$, and the straight line segments joining them give an equivariant homotopy on $|\partial Y_\bullet^\Delta|$. This homotopy is continuous on each product of a closed stratum with a simplex and agrees on faces and degeneracies. As in the proof of Proposition~\ref{prop:BS-duality}, it extends equivariantly over $|Y_\bullet^\Delta|$, since the fixed point space of each stabilizer in $\Cone(T)$ contracts to $0$. Thus $\beta_{\log} \simeq_\Gamma \varphi_{\log} \circ e_{\log}$ as maps of pairs.

We claim that the map induced by $\beta_{\log}$ on homotopy quotients, followed by $\kappa$, is $\beta$. To see this, note that both constructions assign to a chain $\sigma_\bullet$ the diagram $Q(\sigma_\bullet)$ as in Proposition~\ref{prop:flag-map}. The map induced by $e_{\log}$, followed by the projection to $A_g^\Delta$, is $\alpha$, as in Section~\ref{sec:XgDelta}. Under the identifications of relative cohomology used above, the preceding homotopy therefore gives $\alpha^* \circ \varphi_{\log}^* \circ \kappa^* = \beta^*$, as required. \qed

\subsubsection{Proof of Proposition~\ref{prop:equal-boundary-pullbacks}}\label{sec:comparison}

We introduce adapted $\Gamma$-equivariant maps from $\fH$ to $\Cone(T)$ and compare the pullbacks they induce with $\varphi_{\BS}^*$ and $\varphi_{\log}^*$. We prove the existence and homotopy uniqueness of adapted maps and show that each induces both $\varphi_{\BS}^*$ and $\varphi_{\log}^*$ on compactly supported cohomology.

\medskip

For $Z \in \fH$, write $Z = X + iY$. The squared Hodge norm on $\RR^g \oplus \RR^g$ is
\[
q_Z(m, n) := \bigl(m - X n\bigr)^{\mathsf t} Y^{-1}\bigl(m - X n\bigr) + n^{\mathsf t} Yn. 
\]
Note that $q_Z$ varies continuously with $Z$ and satisfies $$q_{\gamma Z}(\gamma v) = q_Z(v)$$ for $\gamma \in \Sp_{2g}(\RR)$. It is identified with the squared Euclidean norm via the symplectic change of coordinates $(m, n) \mapsto \bigl(Y^{-1/2}(m - X n), Y^{1/2} n\bigr)$. Hence, by Cauchy--Schwarz, 
\begin{equation} \label{eq:CS}
|\omega(v, w)| \leq \sqrt{q_Z(v) q_Z(w)}.
\end{equation}
For positive real numbers $s$, let
\[
I_s(Z) = \operatorname{span}_{\QQ} \{v \in \ZZ^{2g} : q_Z(v) < s \}.
\]
By \eqref{eq:CS}, $I_s(Z)$ is spanned by lattice points whose symplectic pairings have absolute value less than $s$. In particular, if $s < 1$ then $I_s(Z)$ is isotropic. 

Fix $0 < \varepsilon < 1$ and set $I(Z) := I_{\varepsilon/2}(Z)$. Let $\fU \subseteq \fH$ be the open $\Gamma$-invariant subset on which $I(Z) \ne 0$.

\begin{definition} 
A continuous $\Gamma$-equivariant map $\psi \colon \fH \to \Cone(T)$ is \emph{adapted} if $\psi(Z) \in T_{I(Z)}$ for $Z \in \fU$.
\end{definition}

\noindent Note that this definition depends on the choice of $\varepsilon$. Since this choice remains fixed throughout, we omit it from the notation and terminology.

We identify $\Cone(T)$ with the set of increasing upper semicontinuous maps $F$ from $[0, 1)$ to the poset of isotropic subspaces of $(\QQ^{2g}, \omega)$. Here upper semicontinuity is equivalent to right continuity: $F(t) = F(s)$ for all $t \geq s$ sufficiently close to $s$. Let $I_0 \subsetneq \cdots \subsetneq I_p$ be the distinct values of $F$. Each fiber is an interval $F^{-1}(I_i) = [a_i, b_i)$, and $F$ corresponds to the point of the simplex indexed by $I_0 \subsetneq \cdots \subsetneq I_p$ with barycentric coordinates $(b_0 - a_0, \ldots, b_p - a_p)$.

Write $F_x$ for the map corresponding to $x \in \Cone(T)$. The cone point corresponds to the constant zero map, and $x \in T$ if and only if $F_x(0) \ne 0$. More generally, $x \in T_I$ if and only if $I \subseteq F_x(0)$.

Suppose $F(s) \subseteq F'(s)$ for every $s \in [0, 1)$. For $t \in [0, 1]$, define
\[
(F *_t F')(s) = \begin{cases}
F(s), & s < t,\\
F'(s), & s \geq t.
\end{cases}
\]
This filtration is increasing and upper semicontinuous. As $t$ runs from $0$ to $1$, it joins $F'$ to $F$.

For arbitrary increasing upper semicontinuous filtrations $F$ and $F'$, define
\[
(F \cap F')(s) = F(s) \cap F'(s).
\]
Their intersection is again increasing and upper semicontinuous and is contained in both filtrations at every $s$. Applying the preceding construction twice therefore joins $F$ to $F \cap F'$ and then to $F'$.

Let $f_0, f_1 \colon S \to \Cone(T)$ be continuous maps from a locally compact Hausdorff space. For $z \in S$, set $F_j = F_{f_j(z)}$, for $j = 0, 1$. Then
\begin{equation}\label{eq:intersection-homotopy}
H(z, t) = \begin{cases}
(F_0 \cap F_1) *_{2t} F_0, & 0 \leq t \leq \tfrac12,\\
(F_0 \cap F_1) *_{2 - 2t} F_1, & \tfrac12 \leq t \leq 1
\end{cases}
\end{equation}
defines a homotopy from $f_0$ to $f_1$ whose value at $(z, \tfrac12)$ is the point corresponding to the filtration $F_0 \cap F_1$. Since $S$ is locally compact, $f_0$ and $f_1$ locally take values in finite subcomplexes. Continuity of $H$ then follows directly from \eqref{eq:intersection-homotopy}. If $f_0$ and $f_1$ are $\Gamma$-equivariant, then so is $H$. Moreover, if $f_0(z)$ and $f_1(z)$ lie in $T_I$, then $H(z, t)$ lies in $T_I$ for every $t$.

\begin{lemma}\label{lem:adapted}
Adapted maps exist, and any two adapted maps are $\Gamma$-equivariantly homotopic through adapted maps.
\end{lemma}

\begin{proof}
For $Z \in \fH$, define
\[
F_Z(s) := \operatorname{span}_{\QQ} \{v \in \ZZ^{2g} : q_Z(v) \leq \tfrac{\varepsilon}{2}(1 + s)\}, \qquad 0 \leq s < 1.
\]
This filtration is increasing and upper semicontinuous, and its values are isotropic by~\eqref{eq:CS}. Let $\psi(Z)$ be the corresponding point in $\Cone(T)$.

It is straightforward to check that $\psi$ is continuous and $\Gamma$-equivariant. Since $I(Z) \subseteq F_Z(0)$, the map $\psi$ is adapted. Finally, the homotopy~\eqref{eq:intersection-homotopy} joins any two adapted maps through adapted maps.
\end{proof}

We now compare the pullback induced by an adapted map with $\varphi_{\BS}^*$ and $\varphi_{\log}^*$.

\begin{lemma}\label{lem:adapted-boundary-comparison}
Let $M$ be $\fH^{\BS}$ or $\fH^{\log}$, and let $\varphi_M \colon (M, \partial M) \to (\Cone(T), T)$ be the corresponding map of pairs. For $y \in \partial M$, let $J_y \subseteq \QQ^{2g}$ be the rational isotropic subspace associated to the stratum containing $y$. Suppose every $y \in \partial M$ has a neighborhood $V_y$ such that $J_y \subseteq I(Z)$ for $Z \in V_y \cap \fH$. Then every adapted map induces $\varphi_M^*$ on compactly supported cohomology.
\end{lemma}

\begin{proof}
Set $V = \fU \cup \partial M$. Since $J_y \ne 0$ and $J_y \subseteq I(Z)$ for $Z \in V_y \cap \fH$, each $V_y$ is contained in $V$. Thus $V$ is open, and $A_g \setminus (\fU/\Gamma)$ is compact because $M/\Gamma$ is compact. An adapted map $\psi$ therefore induces a pullback
\begin{equation}\label{eq:adapted-pullback}
\psi^* \colon H_\Gamma^*(\Cone(T), T) \longrightarrow H_c^*(A_g)
\end{equation}
by pullback to $H_\Gamma^*(\fH, \fU) \cong H^*(A_g, \fU/\Gamma)$, followed by extension of supports.

Choose the neighborhoods $V_y$ equivariantly and small enough that $J_z \subseteq J_y$ for $z \in V_y \cap \partial M$. Then $J_z \subseteq J_y \subseteq I(Z)$ for every $z \in V_y \cap \partial M$ and $Z \in V_y \cap \fH$. Since $\Gamma$ acts properly on $M$ and $\partial M/\Gamma$ is compact, choose a $\Gamma$-equivariant collar $c$ of $\partial M$ sufficiently small that each collar segment lies in one of the $V_y$. It follows that $J_z \subseteq I(c(z, s))$ for every boundary point $z$ and positive collar parameter $s$. Moving the boundary into the interior along this collar and fixing points outside it gives a $\Gamma$-equivariant map $R \colon M \to \fH$. Since the interior of the collar lies in $\fU$, the restriction of $R$ to $\fH$ is homotopic to the identity through $\Gamma$-equivariant maps of $(\fH, \fU)$. Both $\psi R(y)$ and $\varphi_M(y)$ lie in $T_{J_y}$ for $y \in \partial M$, so~\eqref{eq:intersection-homotopy} gives a $\Gamma$-equivariant homotopy between $\psi R$ and $\varphi_M$ through maps of pairs.

Pullback by $\psi R$ factors through $H_\Gamma^*(M, V)$, which equivariant excision identifies with $H_\Gamma^*(\fH, \fU)$, compatibly with the maps to compactly supported cohomology. The homotopy of $R|_{\fH}$ to the identity identifies the induced pullback with~\eqref{eq:adapted-pullback}, and the homotopy between $\psi R$ and $\varphi_M$ identifies it with $\varphi_M^*$, as required.
\end{proof}

\begin{proof}[Proof of Proposition~\ref{prop:equal-boundary-pullbacks}]
By Lemma~\ref{lem:adapted-boundary-comparison}, it suffices to check that each boundary point $y$ in each bordification has a neighborhood $V_y$ such that $J_y \subseteq I(Z)$ for $Z \in V_y \cap \fH$. We first consider $y \in \partial \fH^{\BS}$. 

Suppose $y$ lies in the stratum corresponding to the isotropic flag $J_1 \subsetneq \cdots \subsetneq J_k = J$. Choose symplectic coordinates adapted to the flag and let $Z_0 = i \, \mathrm{Id}_g$. The local product description of the Borel--Serre boundary \cite[\S5.4]{BorelSerre} gives a factorization $Z = u a \ell Z_0$, with $u$ in the unipotent radical of the stabilizer of the flag, $a$ positive diagonal, and $\ell$ in a fixed Levi subgroup of that stabilizer. As $Z \to y$, the factor $u$ converges, and we choose the factor $\ell$ locally continuously so that it also converges. The factor $a$ acts by a scalar $a_j$ on $J_j/J_{j - 1}$, where $J_0 = 0$, with $a_j/a_{j + 1} \to \infty$ for $j < k$ and $a_k \to \infty$. Hence $a^{-1}|_{J_\RR}$ tends to zero. Since $u$ and $\ell$ preserve $J_\RR$ and have bounded inverses, equivariance of the Hodge norm gives $q_Z(v) \to 0$ for $v \in J_\RR$. In particular, each vector in a fixed basis of $J \cap \ZZ^{2g}$ has squared norm less than $\varepsilon/2$ near $y$. Thus $y$ has a neighborhood $V_y$ such that $J \subseteq I(Z)$ for every $Z \in V_y \cap \fH$.

It remains to show that each $y \in \partial \fH^{\log}$ has a neighborhood $V_y$ such that $J_y \subseteq I(Z)$ for $Z \in V_y \cap \fH$. Let $y \in \fH_\sigma^\circ$ and write $J = I_\sigma \otimes \QQ$. As in the Borel--Serre case, it suffices to show that $q_Z(v) \to 0$ for every $v \in J_\RR$ as $Z \in \fH$ approaches $y$.

Choose integral symplectic coordinates in which $J_\RR$ is spanned by the first $r = \dim J$ coordinate vectors in the first summand of $\RR^g \oplus \RR^g$. Write $Z$ in blocks of sizes $r$ and $g - r$. In local toroidal coordinates, $Z_{22}$ tends to a point of $\fH_{g - r}$, $Y_{12}$ remains bounded, and $Y_{11}$ is a sum of the ray generators of $\sigma$ with coefficients tending to infinity, plus a bounded matrix \cite[Chapter~III, \S5]{AMRT}. The ray generators have rational kernels with zero intersection in $J_\RR^\vee$, by the minimality of $I_\sigma$. Their sum is positive definite, so the least eigenvalue of $Y_{11}$ tends to infinity. The restriction of $q_Z$ to $J_\RR$ has matrix $(Y_{11} - Y_{12} Y_{22}^{-1} Y_{21})^{-1}$, which tends to zero. For $r = g$, this matrix is $Y^{-1}$ and again tends to zero. Thus $q_Z(v) \to 0$ for $v \in J_\RR$, as required.
\end{proof}

\subsection{Compatibility of coproducts}\label{sec:block-sums}

For the applications in the next section, we need to know that the coproduct on $H_c^*(\cA)$ respects mixed Hodge structures.

\begin{proposition}\label{prop:block-sums}
The coproduct on $H_c^*(\cA)$ is induced by pullback along the product morphism $m \colon \cA \times \cA \longrightarrow \cA$. Hence, it is a morphism of mixed Hodge structures.
\end{proposition}

\begin{proof}
The assertion is immediate if either factor is $\cA_0$. Fix $g_1, g_2 > 0$ and set $g = g_1 + g_2$. Write $T_n = T_n^\omega(\QQ)$. Under the filtration description, direct sum $\oplus \colon \Cone(T_{g_1}) \times \Cone(T_{g_2}) \to \Cone(T_g)$ is given by $F_{x_1 \oplus x_2}(s) = F_{x_1}(s) \oplus F_{x_2}(s)$.

Let $\widetilde m \colon \fH_{g_1} \times \fH_{g_2} \to \fH_g$ be the block sum map lifting $m$. Using the same $0 < \varepsilon < 1$ in each genus, choose adapted maps $\psi_n \colon \fH_n \to \Cone(T_n)$ for $n \in \{ g_1, g_2, g \}$. Write $\psi_n^*$ for the pullbacks of Lemma~\ref{lem:adapted-boundary-comparison}, and let $\fU_n = \{ Z \in \fH_n : I(Z) \neq 0 \}$. We have an equality of squared Hodge norms $q_{Z_1 \oplus Z_2} = q_{Z_1} \oplus q_{Z_2}$, and hence $I(Z_1 \oplus Z_2) = I(Z_1) \oplus I(Z_2)$. Thus $\widetilde m^{-1}(\fU_g)$ is the union of $\fU_{g_1} \times \fH_{g_2}$ and $\fH_{g_1} \times \fU_{g_2}$. On this union, $\psi_g \widetilde m$ and $\oplus \circ (\psi_{g_1} \times \psi_{g_2})$ take values in $T_{I(Z_1) \oplus I(Z_2)}$, so~\eqref{eq:intersection-homotopy} gives a homotopy through maps of pairs, equivariant for $\Sp_{2g_1}(\ZZ) \times \Sp_{2g_2}(\ZZ)$.

The homotopy gives equality of the induced pullbacks on equivariant relative cohomology. Using that $m$ is proper~\cite[Remark~4.2]{BCGP}, it follows that
\[
(m|_{\cA_{g_1} \times \cA_{g_2}})^* \psi_g^* = (\psi_{g_1}^* \otimes \psi_{g_2}^*) \circ \oplus^*.
\]

The maps $\kappa_n$ commute with direct sum. By Proposition~\ref{prop:E1X}, the induced product on equivariant relative homology corresponds to multiplication on $\cH^\omega$. Proposition~\ref{prop:BS-duality} and Lemma~\ref{lem:adapted-boundary-comparison} identify $\psi_n^* \circ \kappa_n^*$ with the duality isomorphism~\eqref{eq:symplectic-relative-cohomology}. The displayed equality therefore identifies the coproduct on $H_c^*(\cA)$ with $m^*$, as required.
\end{proof}

\section{Betti numbers and mixed Hodge structures}\label{sec:betti-bounds}

In this section, we bound the dimensions of the coinvariant quotient
\[
\cQ := \QQ \otimes_{W_0 H_c^*(\cA)} H_c^*(\cA)
\]
and apply these bounds in small genus. We also give evidence for the conjectures on mixed Hodge structures and multiple zeta values, and prove the exponential-growth statements of Corollary~\ref{cor:betti-growth}.

Write $\cQ_{g, k}$ for the component of $\cQ$ in genus $g$ and cohomological degree $k$. By Proposition~\ref{prop:MM-freeness}, $H_c^*(\cA)$ is isomorphic, as a bigraded module, to $W_0 H_c^*(\cA) \otimes \cQ$. Thus lower bounds for $\dim \cQ_{g, k}$ give lower bounds for the Betti numbers of $\cA_g$.

\begin{proposition}\label{prop:coinvariant-top-weight}
For every $g, k \geq 0$, we have
\[
\dim \cQ_{g, k} \geq \dim \Gr_k^W H_c^k(\cA_g).
\]
\end{proposition}

To prove this bound, we use Deligne's splitting of the weight filtration over $\CC$, which is functorial and compatible with tensor products \cite[Proposition~1.9]{Morgan}. This splitting identifies
\[
\cP := \bigoplus_{g, k \geq 0} \Gr_k^W H_c^k(\cA_g; \CC)
\]
with a subspace of $H_c^*(\cA; \CC)$. By Proposition~\ref{prop:block-sums}, the coproduct respects this splitting. Since $H_c^k(\cA_g; \CC)$ has weights at most $k$, $\cP$ is a connected graded subcoalgebra.

\begin{proof}
The natural map $\cP \to \cQ \otimes_{\QQ} \CC$ is a map of connected graded coalgebras. By Lemma~\ref{lem:coalgebra-injectivity}, it suffices to prove injectivity on primitive elements. By Theorem~\ref{thm:iota-weight-zero} and Proposition~\ref{prop:MM-freeness}, any primitive element in its kernel lies in $W_0 H_c^*(\cA; \CC)$. Every nonzero homogeneous primitive element of $\cP$ has positive weight, so the kernel contains no nonzero primitive element. The map is therefore injective, and comparing dimensions in bidegree $(g, k)$ proves the proposition.
\end{proof}

Propositions~\ref{prop:MM-freeness} and~\ref{prop:coinvariant-top-weight} give lower bounds for Betti numbers, but do not by themselves give the bounds for individual graded pieces of the weight filtration asserted in Proposition~\ref{prop:weight-convolution}. We now prove these bounds using the compatibility of the coproduct with weights.

\begin{proof}[Proof of Proposition~\ref{prop:weight-convolution}]
Work over $\CC$. Restricting the action gives a coalgebra map
\begin{equation}\label{eq:action-pure-weight}
\mu \colon W_0 H_c^*(\cA; \CC) \otimes \cP \longrightarrow H_c^*(\cA; \CC).
\end{equation}
Using Deligne's splitting, project the first factor of the coproduct onto $W_0 H_c^*(\cA; \CC)$ and the second onto $\cP$. The coproduct preserves weight and degree, and the weights of $H_c^k(\cA_g)$ lie between $0$ and $k$, so both projections are coalgebra maps. By symmetry of the product morphisms, the coproduct is graded cocommutative and is therefore also a coalgebra map. We obtain a coalgebra map
\begin{equation} \label{eq:proj-coprod}
H_c^*(\cA; \CC) \xrightarrow{\ \Delta \ } H_c^*(\cA; \CC)^{\otimes 2} \longrightarrow W_0 H_c^*(\cA; \CC) \otimes \cP.
\end{equation}

The action map~\eqref{eq:action-pure-weight}, followed by~\eqref{eq:proj-coprod}, gives a bigraded coalgebra endomorphism of $W_0 H_c^*(\cA; \CC) \otimes \cP$. By Theorem~\ref{thm:iota-weight-zero} and the unit axiom, it fixes $W_0 H_c^*(\cA; \CC) \otimes 1$ and $1 \otimes \cP$, and hence all primitive elements of the tensor product. The composite is injective by Lemma~\ref{lem:coalgebra-injectivity}, so the action \eqref{eq:action-pure-weight} is injective. Each component in fixed genus and cohomological degree is finite dimensional, so the composite is bijective. Thus~\eqref{eq:proj-coprod} is surjective.

Since~\eqref{eq:proj-coprod} preserves weights, its surjectivity gives
\[
\dim \Gr_w^W H_c^m(\cA_g) \geq \sum_{h = 0}^g \dim W_0 H_c^{m - w}(\cA_h) \, \dim \Gr_w^W H_c^w(\cA_{g - h}). \qedhere
\]
\end{proof}

Recall that the tautological ring $R^*(\cA_g) \subseteq H^{2*}(\cA_g; \QQ)$ is generated by the Chern classes $\lambda_i$ of the Hodge bundle. The cycle class map identifies it with the tautological Chow ring \cite[Theorem~13]{GiansiracusaKupersTshishiku}. For $g \geq 1$, it is Gorenstein, with socle
\[
R^{g(g - 1)/2}(\cA_g) = \QQ \cdot \lambda_1 \cdots \lambda_{g - 1}
\]
in cohomological degree $g(g - 1)$ \cite[Theorem~2.5]{vanDerGeerCycles}. Each $R^j(\cA_g)$ is pure of weight $2j$. Write $d_g := g(g + 1)/2$ for the complex dimension of $\cA_g$. For $0 \leq j \leq g(g - 1)/2$, Poincar\'e duality and Proposition~\ref{prop:coinvariant-top-weight} therefore give
\[
\dim \cQ_{g, 2d_g - 2j} \geq \dim R^j(\cA_g).
\]
Since the Poincar\'e polynomial of the tautological ring is $\prod_{i = 1}^{g - 1}(1 + t^{2i})$,
we obtain the coefficientwise bound
\begin{equation}\label{eq:tautological-coinvariant-bound}
\sum_{k \geq 0} \dim \cQ_{g, k} t^k \geq t^{2g} \prod_{i = 1}^{g - 1}(1 + t^{2i}) \qquad (g \geq 1).
\end{equation}
The bigraded module decomposition of Proposition~\ref{prop:MM-freeness} gives
\begin{equation}\label{eq:coinvariant-dimension-convolution}
\dim H_c^k(\cA_g) = \sum_{h = 0}^g \sum_{\ell = 0}^k \dim W_0 H_c^\ell(\cA_h) \, \dim \cQ_{g - h, k - \ell}.
\end{equation}
We now combine these formulas with the known dimensions of weight-zero cohomology and Poincar\'e duality to obtain lower bounds for the Betti numbers of $\cA_g$.

\subsection{Genera at most four}\label{sec:betti-small-genus}

The cohomology of $\cA_g$ for $g \leq 2$ is classical. Hain computed the rational cohomology of $\cA_3$ \cite[Theorem~1]{HainA3}, and its mixed Hodge structure was subsequently determined in \cite[Theorem~1.19]{BCGP}. Writing $E_3$ for the mixed Hodge structure in \eqref{eq:zeta-three-extension}, we have $H_c^6(\cA_3) \cong E_3$ and, by Poincar\'e duality, $H^6(\cA_3) \cong E_3(-3)$. The nonzero cohomology groups for $g \leq 3$, with their mixed Hodge structures, are as follows.
\[
\begin{array}{c|cccc}
g & H^0(\cA_g) & H^2(\cA_g) & H^4(\cA_g) & H^6(\cA_g)\\ \hline
0, 1 & \QQ & 0 & 0 & 0\\
2 & \QQ & \QQ(-1) & 0 & 0\\
3 & \QQ & \QQ(-1) & \QQ(-2) & E_3(-3).
\end{array}
\]
The subspace $W_6 H^6(\cA_3) \cong \QQ(-3)$ is spanned by the tautological class $\lambda_1^3$.

We now apply these bounds in genus four. Write $P(X, t) = \sum_{k \geq 0} b_k(X) t^k$, where $b_k(X)$ is the $k$th rational Betti number. The tautological subring of $H^*(\cA_4)$ has basis
\[
\lambda_1^{\epsilon_1} \lambda_2^{\epsilon_2} \lambda_3^{\epsilon_3}, \qquad \epsilon_i \in \{0, 1\},
\]
and Poincar\'e polynomial $(1 + t^2)(1 + t^4)(1 + t^6)$.

Proposition~\ref{prop:coinvariant-top-weight} gives $\dim \cQ_{1, 2} \geq 1$. Since $W_0 H_c^6(\cA_3) \cong \QQ$, combining \eqref{eq:coinvariant-dimension-convolution} with the tautological bound gives $\dim H_c^8(\cA_4) \geq 2$. Poincar\'e duality and the tautological bounds in the remaining degrees give
\begin{equation}\label{eq:genus-four-betti-bound}
P(\cA_4, t) \geq 1 + t^2 + t^4 + 2t^6 + t^8 + t^{10} + 2t^{12}.
\end{equation}
These bounds agree with those of Hulek and Tommasi \cite[Corollary~3]{HulekTommasiA4}, who also prove $b_{12}(\cA_4) = 2$. Their sum is $9$, equal to the Euler characteristic computed by Ta\"ibi \cite[Proposition~3]{HulekTommasiTopology}.

The product morphism also determines the mixed Hodge structure on $H^{12}(\cA_4)$.

\begin{proposition}\label{prop:genus-four-extension}
The Gysin map for the product morphism $\cA_3 \times \cA_1 \to \cA_4$ induces an isomorphism of mixed Hodge structures
\[
H^6(\cA_3)(-3) \xrightarrow{\ \sim \ } H^{12}(\cA_4).
\]
In particular, $H^{12}(\cA_4)$ is a Tate twist of the $\zeta(3)$ extension.
\end{proposition}

\noindent Equivalently, the coproduct induces an isomorphism of mixed Hodge structures
\[
\Delta_{3, 1} \colon H_c^8(\cA_4) \xrightarrow{\ \sim \ } H_c^6(\cA_3) \otimes H_c^2(\cA_1).
\]
In particular, $H^{12}(\cA_4) \cong E_3(-6)$.

\begin{proof}[Proof of Proposition~\ref{prop:genus-four-extension}]
Let $m \colon \cA_3 \times \cA_1 \to \cA_4$ be the product morphism and $p \colon \cA_3 \times \cA_1 \to \cA_3$ the projection. Consider the Gysin morphism
\[
m_* \colon H^6(\cA_3 \times \cA_1)(-3) \longrightarrow H^{12}(\cA_4).
\]
Its dual, tensored with $\QQ(-10)$, is $\Delta_{3, 1}$. Since $H^*(\cA_1) = \QQ$, pullback by $p$ identifies the source with $H^6(\cA_3)(-3) \cong E_3(-6)$.

The following identity holds in the rational Chow ring of $\cA_4$ \cite[proof of Proposition~37]{COP}:
\[
\lambda_1 m_* 1 = 20 \lambda_1 \lambda_3.
\]
Since $m^* \lambda_1 = p^* \lambda_1$ in cohomology, the projection formula and the relation $\lambda_1^2 = 2 \lambda_2$ give
\[
m_*(p^* \lambda_1^3) = \lambda_1^3 m_* 1 = 20 \lambda_1^3 \lambda_3 = 40 \lambda_1 \lambda_2 \lambda_3 \neq 0.
\]
The displayed basis shows that $m_*$ is nonzero on the weight-$12$ subspace of $E_3(-6)$. Since this is a nonsplit extension of two one-dimensional Hodge structures, $m_*$ is injective. It is an isomorphism because $\dim H^{12}(\cA_4) = 2$.
\end{proof}

\noindent By Poincar\'e duality, Conjecture~\ref{conj:genus-four-mhs} is equivalent to the following description of ordinary cohomology:
\[
H^k(\cA_4) \cong
\begin{cases}
\QQ(-k/2), & k \in \{0, 2, 4, 8, 10\},\\
\QQ(-3)^{ \oplus 2}, & k = 6,\\
E_3(-6), & k = 12,\\
0, & \text{otherwise}.
\end{cases}
\]
It is also equivalent to the vanishing of $H^k(\cA_4)$ for all odd $k$. The predicted compactly supported Euler characteristic agrees, weight by weight, with Ta\"ibi's formula \cite[Theorem~1]{TaibiEuler}.

\subsection{Genus five}\label{sec:betti-genus-five}

By \cite[Table~1]{BCGP}, we have
\[
W_0 H_c^{10}(\cA_5) \cong W_0 H_c^{15}(\cA_5) \cong \QQ.
\]
Proposition~\ref{prop:coinvariant-top-weight} gives $\dim \cQ_{2, 4} \geq 1$ and $\dim \cQ_{2, 6} \geq 1$. Combining these facts with $W_0 H_c^6(\cA_3) \cong \QQ$ and the tautological bound in \eqref{eq:coinvariant-dimension-convolution} gives $\dim H_c^{10}(\cA_5) \geq 3$ and $\dim H_c^{12}(\cA_5) \geq 2$. Poincar\'e duality therefore gives the coefficientwise bound
\begin{equation}\label{eq:genus-five-betti-bound}
P(\cA_5, t) \geq 1 + t^2 + t^4 + 2t^6 + 2t^8 + 2t^{10} + 2t^{12} + 2t^{14} + t^{15} + t^{16} + 2t^{18} + 3t^{20}.
\end{equation}
The alternating sum of these bounds is $19 - 1 = 18$, equal to the Euler characteristic computed by Ta\"ibi \cite[\S7.1]{HulekTommasiTopology}.

The weight-zero classes give
\[
\dim \Gr_{30}^W H^{20}(\cA_5) = \dim \Gr_{30}^W H^{15}(\cA_5) = 1.
\]
Since $H_c^w(\cA_2) \cong \QQ(-w/2)$ for $w = 4, 6$, Proposition~\ref{prop:weight-convolution} and Poincar\'e duality give
\[
\dim \Gr_{26}^W H^{20}(\cA_5) \geq 1 \quad \text{and} \quad \dim \Gr_{24}^W H^{18}(\cA_5) \geq 1.
\]

By Poincar\'e duality, Conjecture~\ref{conj:genus-five-mhs} is equivalent to the following description of ordinary cohomology:
\[
H^k(\cA_5) \cong
\begin{cases}
\QQ(-k/2), & k \in \{0, 2, 4, 16\},\\
\QQ(-k/2)^{ \oplus 2}, & k \in \{6, 8, 10, 12, 14\},\\
\QQ(-15), & k = 15,\\
E_3(-9), & k = 18,\\
E_{3, 5}(-10), & k = 20,\\
0, & \text{otherwise}.
\end{cases}
\]
For every degree and weight, the predicted dimension is a lower bound for the corresponding weight graded piece.

Writing $\LL = [\QQ(-1)]$ in the Grothendieck ring of mixed Hodge structures, the alternating sum of the predicted compactly supported cohomology groups is
\[
\LL^5 \prod_{i = 1}^4(1 + \LL^i) + \LL^2 + \LL^3.
\]
This agrees, weight by weight, with Ta\"ibi's computation of the compactly supported $\ell$-adic Euler characteristic of $\cA_5$ \cite[Theorem~1]{TaibiEuler}.

\subsection{Genus six and higher}\label{sec:betti-higher-genus}

Canning, Oprea, and Pandharipande prove that the algebraic class $[\cA_1 \times \cA_5] \in H^{10}(\cA_6)$ is non-tautological \cite[\S1.6]{COP}. Since the tautological subspace in degree $10$ has dimension three, $\dim \Gr_{10}^W H^{10}(\cA_6) \geq 4$. Poincar\'e duality and Proposition~\ref{prop:coinvariant-top-weight} therefore give
\begin{equation}\label{eq:genus-six-coinvariant-bound}
\sum_{k \geq 0} \dim \cQ_{6, k} t^k \geq t^{12} \prod_{i = 1}^5(1 + t^{2i}) + t^{32}.
\end{equation}
Combining this bound with the weight-zero dimensions in \cite[Table~1]{BCGP} through \eqref{eq:coinvariant-dimension-convolution} gives
\begin{equation}\label{eq:genus-six-betti-bound}
P(\cA_6, t) \geq \prod_{i = 1}^5(1 + t^{2i}) + t^{10} + t^{24} + t^{25} + t^{26} + t^{28} + 3t^{30}.
\end{equation}

Since $H_c^2(\cA_1) \cong \QQ(-1)$ and $\Gr_6^W H_c^6(\cA_3) \cong \QQ(-3)$, Proposition~\ref{prop:weight-convolution}, applied to the weight-zero classes in genera five and three, gives $\dim \Gr_2^W H_c^{12}(\cA_6) \geq 1$ and $\dim \Gr_6^W H_c^{12}(\cA_6) \geq 1$. Together with the tautological class in degree $30$ and $W_0 H_c^{12}(\cA_6) \cong \QQ$, Poincar\'e duality gives
\[
\dim \Gr_w^W H^{30}(\cA_6) \geq 1 \qquad \text{for } w \in \{30, 36, 40, 42\}.
\]

The coefficients on the right side of \eqref{eq:genus-six-betti-bound} sum to $39$ in even degrees and to $1$ in odd degrees, with alternating sum $38$. The Euler characteristic of $\cA_6$ computed by Ta\"ibi is $46$ \cite[\S7.1]{HulekTommasiTopology}. Consequently, the sum of the even-degree Betti numbers is at least $47$, eight more than the lower bound above.

The Euler comparison also gives a lower bound on the coinvariant quotient. The weight-zero Euler characteristics for $g = 0, \ldots, 6$ are $1, 0, 0, 1, 0, 0, 1$ \cite[Table~1]{BCGP}. Taking alternating sums in \eqref{eq:coinvariant-dimension-convolution}, successively through genus six, gives
\[
\sum_{k \geq 0}(-1)^k \dim \cQ_{6, k} = 46 - 5 = 41.
\]
The even-degree components of $\cQ$ in genus six therefore have total dimension at least $41$, eight more than the lower bound of $33$ in \eqref{eq:genus-six-coinvariant-bound}.

The same calculation, using the weight-zero bounds in \cite[Table~1]{BCGP}, gives the following bounds in genera seven through nine:
\begin{align}
\begin{split}
P(\cA_7, t) \geq {} & \prod\nolimits_{i = 1}^6(1 + t^{2i}) + t^{28} + t^{30} + t^{32} + t^{33} + t^{34}\\
& {} + t^{35} + 2t^{36} + 2t^{37} + t^{38} + 2t^{40} + 4t^{42},
\end{split}\label{eq:genus-seven-betti-bound} \displaybreak[1]\\[10pt]
\begin{split}
P(\cA_8, t) \geq {} & \prod\nolimits_{i = 1}^7(1 + t^{2i}) + t^{36} + t^{38} + t^{40} + 3t^{42} + 2t^{44} + t^{45} + 2t^{46}\\
& {} + 2t^{47} + 2t^{48} + t^{49} + 3t^{50} + 4t^{51} + 2t^{52} + 3t^{54} + 6t^{56},
\end{split}\label{eq:genus-eight-betti-bound} \displaybreak[1]\\[10pt]
\begin{split}
P(\cA_9, t) \geq {} & \prod\nolimits_{i = 1}^8(1 + t^{2i}) + t^{42} + t^{44} + t^{45} + t^{46} + 2t^{48} + 3t^{50}\\
& {} + 4t^{52} + 4t^{54} + t^{55} + 4t^{56} + t^{57} + 5t^{58} + 2t^{59} + 4t^{60} + 3t^{61}\\
& {} + 4t^{62} + 3t^{63} + 3t^{64} + 2t^{65} + 5t^{66} + 5t^{67} + 3t^{68} + 4t^{70} + 8t^{72}.
\end{split}\label{eq:genus-nine-betti-bound}
\end{align}
The $t^{32}$ term in \eqref{eq:genus-six-coinvariant-bound} contributes $t^{52}$ to \eqref{eq:genus-nine-betti-bound}, using $W_0 H_c^6(\cA_3) \cong \QQ$.

The right sides of \eqref{eq:genus-seven-betti-bound}--\eqref{eq:genus-nine-betti-bound} have alternating sums $73$, $146$, and $294$, whereas the Euler characteristics of $\cA_7$, $\cA_8$, and $\cA_9$ computed by Ta\"ibi are $104$, $200$, and $528$, respectively \cite[\S7.1]{HulekTommasiTopology}.

\subsection{Multiple zeta values as periods of $\cA_g$}\label{sec:multiple-zeta-values}

In the degree of the socle of the tautological ring, our bounds on the dimensions of the weight graded pieces suggest the presence of a particular subquotient from the theory of mixed Tate motives over $\ZZ$. We explain this first for $H^{56}(\cA_8)$, where the bounds are
\[
\dim \Gr_w^W H^{56}(\cA_8) \geq
\begin{cases}
1, & w \in \{56, 62, 66, 68, 70\},\\
2, & w = 72.
\end{cases}
\]
These dimensions agree with those of a Tate twist of the following mixed Hodge structure associated with multiple zeta values.

Recall that $\cU'_{\leq n}$ is the Hodge realization of the subcomodule of $\cU'$ spanned by words of degree at most $n$. For $n = 8$, these words are
\[
1, \quad f_3, \quad f_5, \quad f_3 f_3, \quad f_7, \quad f_3 f_5, \quad f_5 f_3.
\]
Thus, we have
\[
\Gr^W \cU'_{\leq 8} \cong \QQ \oplus \QQ(-3) \oplus \QQ(-5) \oplus \QQ(-6) \oplus \QQ(-7) \oplus \QQ(-8)^{ \oplus 2}.
\]
After twisting by $\QQ(-28)$, these dimensions agree with the bounds above, as predicted by the conjectured subquotient $\cU'_{\leq 8}(-28)$ of $H^{56}(\cA_8)$. After twisting this subquotient by $\QQ(28)$, we expect the multiple zeta value
\[
\zeta(3, 5) = \sum_{0 < m < n} \frac{1}{m^3 n^5}
\]
to occur as a period.

More generally, Conjecture~\ref{conj:multizeta-socle} predicts that $\cU'_{\leq g}(-g(g - 1)/2)$ is a subquotient of $H^{g(g - 1)}(\cA_g)$. Let $a_n$ be the number of words of total degree $n$ in the letters $f_3, f_5, f_7, \ldots$. The following proposition establishes the corresponding lower bounds for the dimensions of the weight graded pieces.

\begin{proposition}\label{prop:multizeta-weight-bound}
For $g \geq 1$ and $0 \leq n \leq g$, we have
\[
\dim \Gr_{g(g - 1) + 2n}^W H^{g(g - 1)}(\cA_g) \geq a_n.
\]
\end{proposition}

\begin{proof}
By Poincar\'e duality, it is equivalent to prove $\dim \Gr_{2(g - n)}^W H_c^{2g}(\cA_g) \geq a_n$. The classes constructed in \cite{BCGP} generate a free associative subalgebra with one generator in each odd genus at least three, by \cite[Remark~1.3]{Willwacher26}. Counting words in these generators gives $\dim W_0 H_c^{2n}(\cA_n) \geq a_n$. The dual of the socle of the tautological ring in genus $g - n$ contributes a copy of $\QQ(-(g - n))$ to $\Gr_{2(g - n)}^W H_c^{2(g - n)}(\cA_{g - n})$. Taking $m = 2g$ and $w = 2(g - n)$ in Proposition~\ref{prop:weight-convolution}, the summand with $h = n$ gives the required inequality.
\end{proof}

\subsection{Exponential growth or vanishing}\label{sec:betti-exponential-growth}

The numbers $a_n$ in Proposition~\ref{prop:multizeta-weight-bound} also govern the exponential lower bounds on cohomology. Let $S \subseteq W_0 H_c^*(\cA; \RR)$ be the sub-Hopf algebra generated by the primitive classes $\omega^{4i + 1} \in W_0 H_c^{4i + 2}(\cA_{2i + 1}; \RR)$, $i \geq 1$, constructed in \cite{BCGP}. These classes freely generate $S$ as an associative algebra, by \cite[Remark~1.3]{Willwacher26}. By Theorem~\ref{thm:iota-weight-zero}, the map $S \to H_c^*(\cA; \RR)$ given by $a \mapsto a \cdot 1$ is the inclusion. Proposition~\ref{prop:MM-freeness}, applied over $\RR$, therefore shows that $H_c^*(\cA; \RR)$ is a free bigraded $S$-module.

The genus-$g$ part of $S$ lies in cohomological degree $2g$ and has dimension $a_g$. There is one generator in each odd genus at least three, so
\[
\sum_{g \geq 0} a_g z^g = \frac{1 - z^2}{1 - z^2 - z^3}.
\]
The generating function has a unique pole of smallest modulus, which is simple and lies at $1/\beta$, where $\beta = 1.3247 \ldots$ is the real root of $t^3 - t - 1$. Hence $a_g \sim C \beta^g$ for some $C > 0$.

\begin{proof}[Proof of Corollary~\ref{cor:betti-growth}]
Fix an integer $k$. If $H_c^{2g + k}(\cA_g)$ is nonzero for some $g$, then a homogeneous $S$-basis of $H_c^*(\cA; \RR)$ contains an element in bidegree $(h, 2h + k)$ for some $h$. Acting on this element gives
\[
\dim H_c^{2g + k}(\cA_g) \geq a_{g - h} \qquad (g \geq h).
\]
The asymptotic formula for $a_g$ proves the first assertion.

It remains to show that, for each nonnegative $k \neq 1, 3$, the group $H_c^{2h + k}(\cA_h)$ is nonzero for some $h$. If $k$ is even, choose $h \geq 1$ with $h(h - 1) \geq k$. The tautological bound \eqref{eq:tautological-coinvariant-bound} gives $\dim \cQ_{h, 2h + k} \geq 1$, so $H_c^{2h + k}(\cA_h) \neq 0$.

If $k \geq 5$ is odd, choose $h \geq 1$ with $h(h - 1) \geq k - 5$. The same bound gives $\dim \cQ_{h, 2h + k - 5} \geq 1$. Since $W_0 H_c^{15}(\cA_5) \cong \QQ$, the dimension formula \eqref{eq:coinvariant-dimension-convolution} gives
\[
H_c^{2(h + 5) + k}(\cA_{h + 5}) \neq 0.
\]
Poincar\'e duality gives the equivalent statement for ordinary cohomology.
\end{proof}

\begin{remark}\label{rem:symplectic-vanishing}
An affirmative answer to \cite[Question~1.1]{BPS} would imply that $H_c^{2g + k}(\cA_g) = 0$ for all $k < 0$ and all $g$. This is the symplectic analog of the Church--Farb--Putman vanishing conjecture for $\SL_n(\ZZ)$ \cite[Conjecture~2]{CFP}, for which Deng has announced a proof \cite[Theorem~1.1]{DengCFP}. By Corollary~\ref{cor:betti-growth}, any counterexample would force exponential growth of $H^{g(g - 1) - k}(\cA_g)$ for some fixed $k < 0$.

Brown, Panzer, and Portner canonically identify $S$ with the universal enveloping algebra of the motivic Lie algebra of mixed Tate motives over $\ZZ$, after extension of scalars to $\RR$ \cite[Theorem~24]{BrownPanzerPortner26}. The same Lie algebra accounts for the exponential growth of $H^{4g - 6}(\cM_g)$ and the failure of the corresponding vanishing conjecture for mapping class groups \cite[Theorems~1.1 and~2.7]{CGP}.
\end{remark}

\raggedbottom
\bibliographystyle{amsalpha}
\bibliography{ComoduleAlgebra}

\end{document}